\documentclass[a4paper]{amsart}
\usepackage{amsmath,amssymb}

\newtheorem{theorem}{Theorem}[section]
\newtheorem{lemma}[theorem]{Lemma}
\newtheorem{prop}[theorem]{Proposition}
\newtheorem{cor}[theorem]{Corollary}

\theoremstyle{definition}
\newtheorem{definition}[theorem]{Definition}
\newtheorem{example}[theorem]{Example}

\theoremstyle{remark}
\newtheorem{remark}[theorem]{Remark}

\begin{document}
\newcommand{\halfintegers}{\mathbb{Z}+\frac{1}{2}}
\newcommand{\aut}{{\rm aut}}
\newcommand{\modgt}{
\unitlength .1cm
\begin{picture}(5,5)(-5,-5)
\thicklines
\put(0,0){\line(0,1){7}}
\put(0,7){\line(-1,1){4}}
\put(0,7){\line(1,1){4}}
\put(0,0){\line(-1,-1){4}}
\put(0,0){\line(1,-1){4}}
\put(0,0){\circle*{1}}
\put(0,7){\circle*{1}}
\put(1,0){$-a$}
\put(2,6){$a$}
\end{picture}}
\newcommand{\modgtg}{
\unitlength .1cm
\begin{picture}(5,12)(-5,-5)
\thicklines
\put(0,0){\line(0,1){7}}
\put(0,7){\line(-1,1){4}}
\put(0,7){\line(1,1){4}}
\put(0,0){\line(-1,-1){4}}
\put(0,0){\line(1,-1){4}}
\put(0,0){\circle*{1}}
\put(0,7){\circle*{1}}
\put(2,3){$a$}
\put(0,7){\vector(0,-1){5}}
\end{picture}}

\newcommand{\threein}{  
\unitlength .15cm
\begin{picture}(5,12)(-5,-5)
\thicklines
\put(0,0){\line(1,1){5}}
\put(0,0){\line(-1,1){5}}
\put(0,0){\line(0,-1){5}}
\put(0,0){\vector(1,1){4}}
\put(0,0){\vector(-1,1){4}}
\put(0,0){\vector(0,-1){4}}
\put(0,0){\circle*{1}}
\put(4,2){$\mu$}
\put(-5,1.5){$\nu$}
\put(1,-4){$\lambda$}
\put(-2,-1){$v$}
\end{picture}
}
\newcommand{\twoin}{  
\unitlength .15cm
\begin{picture}(5,5)(-5,-5)
\thicklines
\put(0,0){\line(1,1){5}}
\put(0,0){\line(-1,1){5}}
\put(0,0){\line(0,-1){5}}
\put(0,0){\vector(1,1){4}}
\put(0,0){\vector(-1,1){4}}
\put(0,-5){\vector(0,1){3}}
\put(0,0){\circle*{1}}
\put(4,2){$\mu$}
\put(-5,1.5){$\nu$}
\put(1,-4){$\lambda$}
\put(-2,-1){$v$}
\end{picture}
}
\newcommand{\onein}{  
\unitlength .15cm
\begin{picture}(5,5)(-5,-5)
\thicklines
\put(0,0){\line(1,1){5}}
\put(0,0){\line(-1,1){5}}
\put(0,0){\line(0,-1){5}}
\put(0,0){\vector(-1,1){4}}
\put(5,5){\vector(-1,-1){3}}
\put(0,-5){\vector(0,1){3}}
\put(0,0){\circle*{1}}
\put(4,2){$\mu$}
\put(-5,1.5){$\nu$}
\put(1,-4){$\lambda$}
\put(-2,-1){$v$}
\end{picture}
}
\newcommand{\threeout}{  
\unitlength .15cm
\begin{picture}(5,5)(-5,-5)
\thicklines
\put(0,0){\line(1,1){5}}
\put(0,0){\line(-1,1){5}}
\put(0,0){\line(0,-1){5}}
\put(5,5){\vector(-1,-1){3}}
\put(-5,5){\vector(1,-1){3}}
\put(0,-5){\vector(0,1){3}}
\put(0,0){\circle*{1}}
\put(4,2){$\mu$}
\put(-5,1.5){$\nu$}
\put(1,-4){$\lambda$}
\put(-2,-1){$v$}
\end{picture}
}
\newcommand{\treea}{{  
\unitlength .1cm
\begin{picture}(17,17)(-17,-2)
\thicklines
\put(0,0){\line(1,1){15}}
\put(0,0){\line(-1,1){15}}
\put(5,5){\line(-1,1){10}}
\put(0,10){\line(1,1){5}}
\put(0,0){\circle*{1}}
\put(5,5){\circle*{1}}
\put(15,15){\circle*{1}}
\put(-15,15){\circle*{1}}
\put(0,10){\circle*{1}}
\put(-5,15){\circle*{1}}
\put(5,15){\circle*{1}}
\end{picture}
}}
\newcommand{\treeb}{{  
\unitlength .1cm
\begin{picture}(17,12)(-12,-7)
\thicklines
\put(0,0){\line(-1,1){10}}
\put(0,0){\line(1,1){10}}
\put(10,0){\line(1,1){10}}
\put(10,0){\line(-1,1){4}}
\put(0,10){\line(1,-1){4}}
\put(0,0){\circle*{1}}
\put(-10,10){\circle*{1}}
\put(0,10){\circle*{1}}
\put(10,10){\circle*{1}}
\put(20,10){\circle*{1}}
\put(10,0){\circle*{1}}
\end{picture}}}
\newcommand{\treebb}{{  
\unitlength .1cm
\begin{picture}(17,12)(-12,-7)
\thicklines
\put(0,0){\line(-1,1){10}}
\put(0,0){\line(1,1){10}}
\put(10,0){\line(1,1){10}}
\put(10,0){\line(-1,1){4}}
\put(0,10){\line(1,-1){4}}
\put(0,0){\circle{1.5}}
\put(-10,10){\circle*{1}}
\put(0,10){\circle*{1}}
\put(10,10){\circle*{1}}
\put(20,10){\circle*{1}}
\put(10,0){\circle{1.5}}
\end{picture}}}

\newcommand{\treec}{{  
\unitlength .1cm
\begin{picture}(17,17)(-17,-2)
\thicklines
\put(0,0){\line(1,1){15}}
\put(0,0){\line(-1,1){15}}
\put(0,10){\line(-1,1){5}}
\put(0,10){\line(1,1){5}}
\put(0,0){\circle*{1}}
\put(-15,15){\circle*{1}}
\put(-5,15){\circle*{1}}
\put(5,15){\circle*{1}}
\put(15,15){\circle*{1}}
\put(0,10){\circle*{1}}
\end{picture}}}
\newcommand{\treecc}{{  
\unitlength .1cm
\begin{picture}(17,17)(-17,-2)
\thicklines
\put(0,0){\line(1,1){15}}
\put(0,0){\line(-1,1){15}}
\put(0,10){\line(-1,1){5}}
\put(0,10){\line(1,1){5}}
\put(0,0){\circle{1.5}}
\put(-15,15){\circle*{1}}
\put(-5,15){\circle*{1}}
\put(5,15){\circle*{1}}
\put(15,15){\circle*{1}}
\put(0,10){\circle{1.5}}
\end{picture}}}

\newcommand{\treed}{{  
\unitlength .1cm
\begin{picture}(12,7)(-2,-2)
\thicklines
\put(5,0){\line(-1,1){5}}
\put(5,0){\line(1,1){5}}
\put(0,5){\circle*{1}}
\put(10,5){\circle*{1}}
\put(5,0){\circle*{1}}
\end{picture}
}}
\newcommand{\treedd}{{  
\unitlength .1cm
\begin{picture}(12,7)(-2,-2)
\thicklines
\put(5,0){\line(-1,1){5}}
\put(5,0){\line(1,1){5}}
\put(0,5){\circle*{1}}
\put(10,5){\circle*{1}}
\put(5,0){\circle{1.5}}
\end{picture}
}}
\newcommand{\treee}{{
\unitlength .1cm
\begin{picture}(12,12)(-14,-2)
\thicklines
\put(0,0){\line(1,1){10}}
\put(0,0){\line(-1,1){10}}
\put(5,5){\line(-1,1){5}}
\put(-10,10){\circle*{1}}
\put(0,10){\circle*{1}}
\put(10,10){\circle*{1}}
\put(5,5){\circle*{1}}
\put(0,0){\circle*{1}}
\end{picture}
}}
\newcommand{\treeee}{{
\unitlength .1cm
\begin{picture}(12,12)(-14,-2)
\thicklines
\put(0,0){\line(1,1){10}}
\put(0,0){\line(-1,1){10}}
\put(5,5){\line(-1,1){5}}
\put(-10,10){\circle*{1}}
\put(0,10){\circle*{1}}
\put(10,10){\circle*{1}}
\put(5,5){\circle*{1}}
\put(0,0){\circle{1.5}}
\end{picture}
}}
\newcommand{\treef}{{
\unitlength .1cm
\begin{picture}(12,12)(-14,-2)
\thicklines
\put(0,0){\line(1,1){10}}
\put(0,0){\line(-1,1){10}}
\put(-10,10){\circle*{1}}
\put(0,10){\circle*{1}}
\put(10,10){\circle*{1}}
\put(0,0){\circle*{1}}
\end{picture}
}}
\newcommand{\treeff}{{
\unitlength .1cm
\begin{picture}(12,12)(-14,-2)
\thicklines
\put(0,0){\line(1,1){10}}
\put(0,0){\line(-1,1){10}}
\put(-10,10){\circle*{1}}
\put(0,10){\circle*{1}}
\put(10,10){\circle*{1}}
\put(0,0){\circle{1.5}}
\end{picture}
}}
\newcommand{\ancomb}[2]{{
\unitlength .1cm
\begin{picture}(70,45)(-50,-15)
\put(-43,-12){#1}
\put(7,-12){#2}
\put(-8,18){\treed}
\put(45,-7){\treee}
\put(24,-15){$e_1$}
\put(-26,-15){$e_3$}
\put(-1,16){$e_2$}
\put(60,-11){$e_4$}
\qbezier(30,5)(45,20)(60,5)
\qbezier(-10,5)(0,15)(10,5)
\qbezier(-5,25)(15,40)(40,5)
\qbezier(5,25)(-15,40)(-40,5)
\put(-24,19){$\mu^2_1$}
\put(21,19){$\mu^1_1$}
\put(-1,5){$\mu^3_1$}
\put(45,14){$\nu^{1,1}_1$}
\end{picture}
}}

\title[Integrality of Gopakumar--Vafa Invariants]
{Integrality of 
Gopakumar--Vafa Invariants of Toric Calabi--Yau Threefolds}
\author{Yukiko Konishi}
\thanks{Communicated by K. Saito. Received April 15, 2005.
Revised October 24, 2005.}
\thanks{{2000} {\it Mathematics Subject Classification} 
Primary 14N35; Secondary 05E05.}
\thanks{Research Institute for Mathematical Sciences,
Kyoto University, Kyoto 606-8502, Japan.}
\thanks{E-mail: konishi@kurims.kyoto-u.ac.jp}
\begin{abstract}
The Gopakumar--Vafa invariants are
numbers defined as certain linear combinations of 
the Gromov--Witten invariants.
We prove that
the GV invariants
of a toric Calabi--Yau threefold
are integers and 
that the invariants for high genera vanish.
The proof of the integrality 
is based on elementary number theory and
that of 
the vanishing uses the 
operator formalism and the exponential formula.
\end{abstract}

\maketitle
\section{Introduction}%
A toric Calabi--Yau (TCY)
threefold is
a three-dimensional smooth toric variety of finite type,
whose canonical bundle is trivial.
For example,
the total space of 
the rank two vector bundle over $\mathbb{P}^1$,
${\mathcal O}(a_1)\oplus {\mathcal O}(a_2)\to \mathbb{P}^1$,
such that $a_1+a_2=-2$
and
the total space of the canonical bundle
of a smooth toric surface
are TCY threefolds.
%
%

Thanks to the duality 
of open and closed strings,
a procedure 
to write down
the partition function 
of the 0-pointed Gromov--Witten (GW) invariants
of any TCY threefold $X$
became available \cite{akmv}.
By the partition function, we mean the exponential 
of the the generating function.
One only has to draw a labeled planar graph from
the fan of $X$
and combine a certain quantity
according to the shape and the labels of the graph.
See 
\cite{z1}\cite{llz1}\cite{llz2}\cite{lllz} for 
the mathematical formulation and the proof.
In this article, we call the graph the {toric graph}
of $X$ and
refer to
the quantity  as the { three point function}.

One open problem concerning the Calabi--Yau threefold
is the Gopakumar--Vafa (GV) conjecture \cite{gv}.
We define the Gopakumar--Vafa invariants
as certain linear combinations
of the GW invariants in the manner of \cite{bp}.
One statement of 
the conjecture is that
the GV invariants
are integers and  
that only finite number of them
are nonzero (in a given homology class).
This is remarkable given that the GW invariants themselves
are, in general, not integers but rational numbers.
Other statement is that 
the GV invariants are equal to ``the number of BPS states''
in the M-theory compactified on the TCY threefold.
A mathematical formulation in this direction
was proposed in \cite{hst}.
Recently,
the studies
using the relation to the instanton counting
appeared
\cite{lilz}\cite{ak}.

The first statement of the GV conjecture
was proved  by Peng \cite{peng} 
in the case of the canonical bundles of Fano toric surfaces.
The aim of this article
is to prove it
for general TCY threefolds.
We put the problem in a combinatorial setting
and %
prove the combinatorial version of the statement.
The proof consists of two parts
%
corresponding to
the integrality and the vanishing for high genera.
The proof of the former is based on
elementary number theory
and basically the same as
that of \cite{peng}.
The proof of the latter uses the operator formalism and
the exponential formula.
It is the generalization of the results of \cite{ko}.

The organization of the paper is as follows.
In section \ref{PF},
we define a generalization of the toric graph,
the partition function and the free energy.
In section \ref{mainresults},
we state main results.
In section \ref{TCY},
we explain
that the first statement of the GV conjecture
follows from these results.
In sections \ref{PI} and \ref{PV},
we give proofs of the integrality and the vanishing, respectively.
Appendix contains a proof of a lemma.

\begin{center}{\bf Acknowledgement}\end{center} %
The author thanks H.Kanno, A.Kato, H. Awata, A. Takahashi,
A. Tsuchiya
and S. Hosono for valuable comments.

\section{Partition Function}\label{PF}
In this section, we first define
the notion of the generalized toric (GT) graph.
Then we
introduce the three point function
and
define the partition function and the free energy of the GT graph.
%

\subsection{Generalized Toric Graph}
Throughout this article,
we assume that
a graph has the finite edge set and vertex set 
and has no self-loop.

A flag $f$ is a pair of a vertex $v$ and an edge $e$
such that $e$ is incident on $v$.
The flag whose edge is the same as $f$ and
vertex is the other endpoint of the edge is denoted by $-f$.
\newcommand{\flag}{{  
\unitlength .15cm
\begin{picture}(17,5)(-2,-3)
\thicklines
\put(0,0){\line(1,0){10}}
\put(0,0){\circle*{1}}
\put(10,0){\circle*{1}}
\put(-1,1){$v$}
\put(10,1){$v'$}
\put(4,1){$e$}
\put(17,2){$f=(v,e)$}
\put(15,-2){$-f=(v',e)$}
\end{picture}
}}
\begin{equation}\notag
\flag
\end{equation}

A connected planar graph $\Gamma$ 
is  a {\em trivalent planar graph} if
all vertices are either trivalent
or univalent.
The set of trivalent vertices  is denoted by
$V_3(\Gamma)$.
The set of edges whose two endpoints are both trivalent 
is denoted by $E_3(\Gamma)$.
The set of flags whose edges are in $E_3(\Gamma)$
is denoted by $F_3(\Gamma)$.
%
\begin{definition}
A trivalent planar graph
with a label $n_f\in\mathbb{Z}$ 
on every flag $f\in F_3(\Gamma)$
%
together with
a drawing into $\mathbb{R}^2$
is a {\em generalized toric graph} (GT graph)
if it satisfies the following conditions.
\begin{enumerate}
\item
$n_f=-n_{-f}$.
\item The drawing has no crossing.
\end{enumerate}
$n_f$ is called the {\em framing} of the flag $f$.
\end{definition}
%

Since $n_f=-n_{-f}$,
assigning framings is the same as assigning
each edge an integer and a direction.
Therefore, 
we add an auxiliary direction to every edge $e\in E_3(\Gamma)$
and redraw the graph 
as follows.

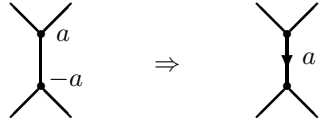
\begin{equation}\notag
\raisebox{-.7cm}{\modgt}
\qquad\qquad
\Rightarrow
\qquad
\raisebox{-.7cm}{\modgtg}
\end{equation}
%
The direction of the edge is taken arbitrarily.
The label on an edge $e$ is denoted by $n_e$.

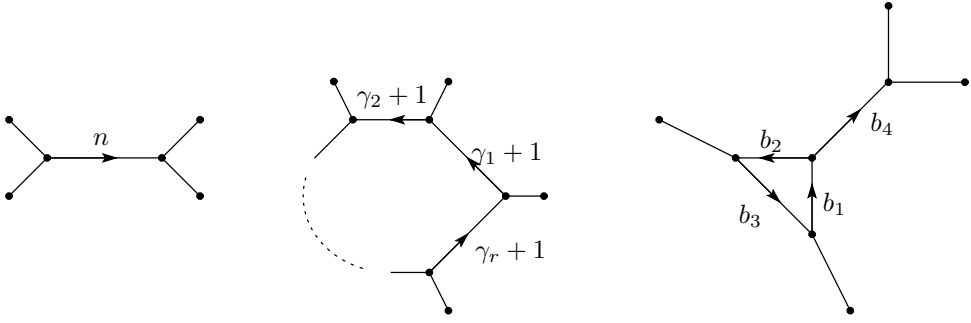
\begin{figure}
\begin{center}
\unitlength 0.1in
\begin{picture}( 50.0600, 16.0300)( 12.0000,-22.0300)
%
\special{pn 8}%
\special{pa 1200 1200}%
\special{pa 1400 1400}%
\special{fp}%
\special{pa 1400 1400}%
\special{pa 1200 1600}%
\special{fp}%
\special{pa 1400 1400}%
\special{pa 2000 1400}%
\special{fp}%
\special{pa 2000 1400}%
\special{pa 2200 1200}%
\special{fp}%
\special{pa 2000 1400}%
\special{pa 2200 1600}%
\special{fp}%
%
\special{pn 8}%
\special{pa 1400 1400}%
\special{pa 1760 1400}%
\special{fp}%
\special{sh 1}%
\special{pa 1760 1400}%
\special{pa 1694 1380}%
\special{pa 1708 1400}%
\special{pa 1694 1420}%
\special{pa 1760 1400}%
\special{fp}%
\put(16.4000,-13.3000){\makebox(0,0)[lb]{$n$}}%
%
\special{pn 20}%
\special{sh 1}%
\special{ar 1400 1400 10 10 0  6.28318530717959E+0000}%
\special{sh 1}%
\special{ar 1200 1200 10 10 0  6.28318530717959E+0000}%
\special{sh 1}%
\special{ar 1200 1600 10 10 0  6.28318530717959E+0000}%
\special{sh 1}%
\special{ar 2200 1600 10 10 0  6.28318530717959E+0000}%
\special{sh 1}%
\special{ar 2200 1200 10 10 0  6.28318530717959E+0000}%
\special{sh 1}%
\special{ar 2000 1400 10 10 0  6.28318530717959E+0000}%
%
\special{pn 8}%
\special{pa 3200 1200}%
\special{pa 3400 1200}%
\special{fp}%
\special{pa 3400 1200}%
\special{pa 3800 1600}%
\special{fp}%
\special{pa 3800 1600}%
\special{pa 3400 2000}%
\special{fp}%
\special{pa 3400 2000}%
\special{pa 3200 2000}%
\special{fp}%
%
\special{pn 8}%
\special{pa 3210 1200}%
\special{pa 3000 1200}%
\special{fp}%
\special{pa 3000 1200}%
\special{pa 2800 1400}%
\special{fp}%
%
\special{pn 8}%
\special{ar 3200 1600 460 400  1.9028558 1.9307628}%
\special{ar 3200 1600 460 400  2.0144837 2.0423907}%
\special{ar 3200 1600 460 400  2.1261116 2.1540186}%
\special{ar 3200 1600 460 400  2.2377395 2.2656465}%
\special{ar 3200 1600 460 400  2.3493674 2.3772744}%
\special{ar 3200 1600 460 400  2.4609953 2.4889023}%
\special{ar 3200 1600 460 400  2.5726232 2.6005302}%
\special{ar 3200 1600 460 400  2.6842511 2.7121581}%
\special{ar 3200 1600 460 400  2.7958791 2.8237860}%
\special{ar 3200 1600 460 400  2.9075070 2.9354139}%
\special{ar 3200 1600 460 400  3.0191349 3.0470418}%
\special{ar 3200 1600 460 400  3.1307628 3.1586697}%
\special{ar 3200 1600 460 400  3.2423907 3.2702977}%
\special{ar 3200 1600 460 400  3.3540186 3.3819256}%
%
\special{pn 8}%
\special{pa 3800 1600}%
\special{pa 4010 1600}%
\special{fp}%
%
\special{pn 8}%
\special{pa 3400 2000}%
\special{pa 3500 2200}%
\special{fp}%
%
\special{pn 8}%
\special{pa 3500 1000}%
\special{pa 3400 1200}%
\special{fp}%
%
\special{pn 8}%
\special{pa 3000 1200}%
\special{pa 2900 1000}%
\special{fp}%
%
\special{pn 20}%
\special{sh 1}%
\special{ar 3000 1200 10 10 0  6.28318530717959E+0000}%
\special{sh 1}%
\special{ar 3400 1200 10 10 0  6.28318530717959E+0000}%
\special{sh 1}%
\special{ar 3800 1600 10 10 0  6.28318530717959E+0000}%
\special{sh 1}%
\special{ar 3400 2000 10 10 0  6.28318530717959E+0000}%
\special{sh 1}%
\special{ar 3500 2200 10 10 0  6.28318530717959E+0000}%
\special{sh 1}%
\special{ar 4000 1600 10 10 0  6.28318530717959E+0000}%
\special{sh 1}%
\special{ar 3500 1000 10 10 0  6.28318530717959E+0000}%
\special{sh 1}%
\special{ar 2900 1000 10 10 0  6.28318530717959E+0000}%
%
\special{pn 8}%
\special{pa 3400 2000}%
\special{pa 3600 1800}%
\special{fp}%
\special{sh 1}%
\special{pa 3600 1800}%
\special{pa 3540 1834}%
\special{pa 3562 1838}%
\special{pa 3568 1862}%
\special{pa 3600 1800}%
\special{fp}%
%
\special{pn 8}%
\special{pa 3800 1600}%
\special{pa 3600 1400}%
\special{fp}%
\special{sh 1}%
\special{pa 3600 1400}%
\special{pa 3634 1462}%
\special{pa 3638 1438}%
\special{pa 3662 1434}%
\special{pa 3600 1400}%
\special{fp}%
%
\special{pn 8}%
\special{pa 3400 1200}%
\special{pa 3200 1200}%
\special{fp}%
\special{sh 1}%
\special{pa 3200 1200}%
\special{pa 3268 1220}%
\special{pa 3254 1200}%
\special{pa 3268 1180}%
\special{pa 3200 1200}%
\special{fp}%
\put(36.4000,-19.5000){\makebox(0,0)[lb]{$\gamma_r+1$}}%
\put(36.2000,-14.4000){\makebox(0,0)[lb]{$\gamma_1+1$}}%
\put(30.2000,-11.4000){\makebox(0,0)[lb]{$\gamma_2+1$}}%
%
\special{pn 8}%
\special{pa 5004 1400}%
\special{pa 5404 1400}%
\special{pa 5404 1800}%
\special{pa 5404 1800}%
\special{pa 5004 1400}%
\special{fp}%
%
\special{pn 8}%
\special{pa 5404 1800}%
\special{pa 5604 2200}%
\special{fp}%
%
\special{pn 8}%
\special{pa 4604 1200}%
\special{pa 5004 1400}%
\special{fp}%
%
\special{pn 8}%
\special{pa 5804 604}%
\special{pa 5804 1004}%
\special{fp}%
\special{pa 5804 1004}%
\special{pa 6204 1004}%
\special{fp}%
%
\special{pn 20}%
\special{sh 1}%
\special{ar 5804 604 10 10 0  6.28318530717959E+0000}%
%
\special{pn 20}%
\special{sh 1}%
\special{ar 4604 1200 10 10 0  6.28318530717959E+0000}%
%
\special{pn 20}%
\special{sh 1}%
\special{ar 5604 2200 10 10 0  6.28318530717959E+0000}%
%
\special{pn 20}%
\special{sh 1}%
\special{ar 5404 1800 10 10 0  6.28318530717959E+0000}%
%
\special{pn 20}%
\special{sh 1}%
\special{ar 5004 1400 10 10 0  6.28318530717959E+0000}%
%
\special{pn 20}%
\special{sh 1}%
\special{ar 5404 1400 10 10 0  6.28318530717959E+0000}%
%
\special{pn 20}%
\special{sh 1}%
\special{ar 5804 1004 10 10 0  6.28318530717959E+0000}%
%
\special{pn 20}%
\special{sh 1}%
\special{ar 6204 1004 10 10 0  6.28318530717959E+0000}%
%
\special{pn 8}%
\special{pa 5404 1800}%
\special{pa 5404 1540}%
\special{fp}%
\special{sh 1}%
\special{pa 5404 1540}%
\special{pa 5384 1608}%
\special{pa 5404 1594}%
\special{pa 5424 1608}%
\special{pa 5404 1540}%
\special{fp}%
%
\special{pn 8}%
\special{pa 5404 1400}%
\special{pa 5134 1400}%
\special{fp}%
\special{sh 1}%
\special{pa 5134 1400}%
\special{pa 5200 1420}%
\special{pa 5186 1400}%
\special{pa 5200 1380}%
\special{pa 5134 1400}%
\special{fp}%
%
\special{pn 8}%
\special{pa 5004 1400}%
\special{pa 5234 1630}%
\special{fp}%
\special{sh 1}%
\special{pa 5234 1630}%
\special{pa 5200 1570}%
\special{pa 5196 1592}%
\special{pa 5172 1598}%
\special{pa 5234 1630}%
\special{fp}%
%
\special{pn 8}%
\special{pa 5804 1000}%
\special{pa 5404 1400}%
\special{fp}%
%
\special{pn 8}%
\special{pa 5404 1400}%
\special{pa 5654 1150}%
\special{fp}%
\special{sh 1}%
\special{pa 5654 1150}%
\special{pa 5592 1184}%
\special{pa 5616 1188}%
\special{pa 5620 1212}%
\special{pa 5654 1150}%
\special{fp}%
\put(54.6000,-17.0000){\makebox(0,0)[lb]{$b_1$}}%
\put(51.3000,-13.7000){\makebox(0,0)[lb]{$b_2$}}%
\put(50.2000,-17.6000){\makebox(0,0)[lb]{$b_3$}}%
\put(57.1000,-12.8000){\makebox(0,0)[lb]{$b_4$}}%
\end{picture}%
\end{center}
\caption{
Examples of the GT graph
($n,\gamma_1,\ldots,\gamma_r,b_1,\ldots,b_4\in\mathbb{Z}$).
}
\label{exgts}
\end{figure}

Examples of the GT graphs are shown in figure \ref{exgts}.

%


\subsection{Partition and Notations}
We summarize notations (mainly) on partitions
(\cite{mac}).
%

A {\em partition} is  a non-increasing sequence 
$\lambda=(\lambda_1,\lambda_2,\ldots)$ of 
nonnegative integers containing only finitely many nonzero
terms.
The nonzero $\lambda_i$'s  are called the parts.
The number of parts is the {\em length} of $\lambda$, 
denoted by $l(\lambda)$. 
The sum of the parts is the {\em weight} of 
$\lambda$, denoted by $|\lambda|$: $|\lambda|=\sum_i\lambda_i$.
If $|\lambda|=d$, $\lambda$ is a partition of $d$. 
The set of all partitions of $d$ is denoted by ${\mathcal P}_d$
and the set of all partitions by ${\mathcal P}$. 
Let $m_k(\lambda)=\#\{\lambda_i:\lambda_i=k\}$ 
be the {\em multiplicity}
of $k$ where
$\#$ denotes the number of elements of a finite set. 
Let $\aut(\lambda)$ be the symmetric group 
acting as the permutations of  the equal parts of $\lambda$:
$\aut(\lambda)\cong \prod_{k\geq 1}\mathfrak{S}_{m_k(\lambda)}$.
Then $\#\aut(\lambda)=\prod_{k\geq 1}m_k(\lambda)!$.
We define
\begin{equation}\notag
z_{\lambda}=\prod_{i=1}^{l(\lambda)}\lambda_i \cdot\#\aut(\lambda),
\end{equation}
which is the number of the centralizers of the conjugacy class 
associated to $\lambda$.

A partition $\lambda=(\lambda_1,\lambda_2,\ldots)$ is identified as 
the Young diagram with $\lambda_i$ boxes in the $i$-th row 
$(1\leq i\leq l(\lambda))$.
The Young diagram with $\lambda_i$ boxes in the $i$-th column is 
its {\em transposed} Young diagram. The corresponding partition
is called the {\em conjugate partition} and denoted by $\lambda^t$.
Note that
$\lambda^t_i=\sum_{k\geq i}m_k(\lambda)$.

We define
\begin{equation}\notag
\kappa(\lambda)=\sum_{i=1}^{l(\lambda)}\lambda_i(\lambda_i-2i+1).
\end{equation}
This is equal to twice the sum of contents $\sum_{x\in \lambda}c(x)$ 
where $c(x)=j-i$ for the box $x$ at the $(i,j)$-th place 
in the Young diagram $\lambda$.
Thus, $\kappa(\lambda)$ is always even and 
satisfies $\kappa(\lambda^t)=-\kappa(\lambda)$.

$\mu\cup\nu$ denotes the partition whose parts 
are $\mu_1,\ldots,\mu_{l(\mu)}, 
\nu_1,\ldots,\nu_{l(\mu)}$ and 
$k\mu$ the partition $(k\mu_1,k\mu_2,\ldots)$ 
for $k\in {\mathbb{N}}$.

For a finite set of integers $s=(s_1,s_2,\ldots,s_l)$,
we use the following notations.
\begin{equation}\notag
|s|=\sum_i s_i.
\end{equation}
%
%
When $s$ has at least one nonzero element, we define
\begin{equation}\notag
{\rm gcd}(s)=\text{the greatest common divisor of }\{|s_i|,s_i\neq 0\}
\end{equation}
where $|s_i|$ is the absolute value of $s_i$.

Throughout this paper, we use the letter $q$ for a variable.
We define
\begin{equation}\notag
[k]=q^{\frac{k}{2}}-q^{-\frac{k}{2}} \qquad (k\in\mathbb{Q}),
\end{equation}
which is called the {\em $q$-number}.
For a partition $\lambda$ and a finite set $s$ as 
above, we use the shorthand notations
\begin{equation}\notag
[\lambda]=\prod_{i=1}^{l(\lambda)}[\lambda_i],
\qquad
[s]=[s_1]\ldots[s_l].
\end{equation}
%

\subsection{Three Point Function}
%
Let $q^{\rho}$ and $q^{\lambda+\rho}$ be the following
(infinite) sequences:
\begin{equation}\notag
q^{\rho}=(q^{-i+\frac{1}{2}})_{i\geq 1},\qquad
q^{\lambda+\rho}=(q^{\lambda_i-i+\frac{1}{2}})_{i\geq 1}.
\end{equation}
The Schur function 
and skew-Schur function are denoted by
$s_{\lambda}$ and $s_{\lambda/\mu}$.

\begin{definition}
Let $(\lambda^1,\lambda^2,\lambda^3)$
be a triple of partitions.
The {\em three point function} is
\begin{equation}\notag
C_{\lambda^1,\lambda^2,\lambda^3}(q)=
q^{\frac{\kappa(\lambda^3)}{2}}
s_{\lambda^2}(q^{\rho})\sum_{\eta\in{\mathcal P}}
s_{\lambda^1/\eta}(q^{\lambda^{2t}+\rho})
s_{\lambda^{3t}/\eta}(q^{\lambda^2+\rho}).
\end{equation}
%
\end{definition}%
This is a rational function in $q^{\frac{1}{2}}$.
An important property of the three point function is 
the cyclic symmetry: 
\begin{equation}\notag
C_{\lambda^1,\lambda^2,\lambda^3}(q)
=C_{\lambda^2,\lambda^3,\lambda^1}(q)
=C_{\lambda^3,\lambda^1,\lambda^2}(q).
\end{equation}
See \cite{orv} for a
proof.
Various identities can be found in \cite{z2}.
Since the variables $q^{\rho}$ and $q^{\lambda+\rho}$ 
are infinite sequences,
let us explain how to compute the 
(skew-) Schur function.
For a sequence of variables  $x=(x_1,x_2,\ldots)$,
the elementary symmetric function
$e_i(x)$ ($i\geq 0$)
and 
the completely symmetric function $h_i(x)$ ($i\geq 0$)
are obtained from
the generating functions:
\begin{equation}\notag
\sum_{i=0}^{\infty}e_i(x)z^i=
\prod_{i=0}^{\infty}(1+x_i z),
\qquad
\sum_{i=0}^{\infty}h_i(x)z^i=
\prod_{i=0}^{\infty}(1-x_i z)^{-1}.
\end{equation}
The skew-Schur function $s_{\mu/\nu}(x)$
%
is written in terms of 
$e_i(x)$ or $h_i(x)$:
\begin{equation}\label{she}
s_{\mu/\nu}(x)=
\det
\big(e_{\mu_i^t-\nu_j^t-i+j}(x)
\big)_{1\leq i,j\leq l(\mu^t)}
=
\det
\big(h_{\mu_i-\nu_j-i+j}(x)
\big)_{1\leq i,j\leq l(\mu)}.
\end{equation}
In the determinants, 
$h_{-i}(x)$ and $e_{-i}(x)$ ($i>0$) are assumed to be zero.
For variables $q^{\rho}$,
we can compute
the elementary and the completely symmetric functions
by using
the identities \cite{mac}:
\begin{equation}\notag
\begin{split}
&
\prod_{i=0}^{\infty}(1+q^i z)=
\sum_{i=0}^{\infty}\frac{q^{\frac{i(i-1)}{2}}}
{(1-q)\ldots (1-q^i)}
z^i,
\\
&
\prod_{i=0}^{\infty}(1-q^i z)^{-1}=
\sum_{i=0}^{\infty}\frac{1}
{(1-q)\ldots(1-q^i)}
\,z^i.
\end{split}
\end{equation}
Therefore, for the variable $q^{\rho}$,
\begin{equation}\label{eh1}
e_i(q^{\rho})=\frac{q^{-\frac{i(i-1)}{4}}}{[1]\ldots[i]},
\qquad
h_i(q^{\rho})=\frac{q^{\frac{i(i-1)}{4}}}{[1]\ldots[i]}.
\end{equation}
For the variable $q^{\lambda+\rho}$,
$e_i(q^{\lambda+\rho})$ and $h_i(q^{\lambda+\rho})$
are computed from the generating functions:
\begin{equation}\label{eh2}
\begin{split}
\sum_{i=0}^{\infty}e_i(q^{\lambda+\rho})z^i
&=\prod_{i=1}^{l(\lambda)}
\frac{1+q^{\lambda_i-i+\frac{1}{2}}z}{1+q^{-i+\frac{1}{2}}z}
\cdot
\Big(\sum_{k= 0}^{\infty}e_k(q^{\rho})z^k \Big),
\\
\sum_{i=0}^{\infty}h_i(q^{\lambda+\rho})z^i
&=\prod_{i=1}^{l(\lambda)}
\frac{1-q^{-i+\frac{1}{2}}z}{1-q^{\lambda_i-i+\frac{1}{2}}z}
\cdot
\Big(\sum_{k=0}^{\infty}h_k(q^{\rho})z^k \Big).
\end{split}
\end{equation}
In this way, we can explicitly compute the skew-Schur functions
and the three point functions.
Here are some examples of three point functions.
\begin{equation}\notag
\begin{split}
&C_{(1),\emptyset,\emptyset}(q)=\frac{1}{[1]},
\quad
C_{(2),\emptyset,\emptyset}(q)=\frac{q^2}{(q-1)(q^2-1)},
\quad
C_{(1,1),\emptyset,\emptyset}(q)=\frac{q}{(q-1)(q^2-1)},
\\
&C_{(1),(1),\emptyset}(q)=\frac{q^2-q+1}{(1-q)^2},
\quad
C_{(1),(1),(1)}(q)=\frac{q^4-q^3+q^2-q+1}{q^{\frac{1}{2}}(q-1)^3}.
\end{split}
\end{equation}
More examples can be found in \cite{akmv}, section 8.

\subsection{Partition Function}
First we set some notations.
Consider a GT graph $\Gamma$.
\begin{itemize}
\item
We associate one formal variable to every edge $e\in E_3(\Gamma)$.
The variable associated to $e$ is denoted by $Q_e$.
$\Vec{Q}=(Q_e)_{e\in E_3(\Gamma)}$.
\item
A {\em degree} is 
a set $\Vec{d}=(d_e)_{e\in E_3(\Gamma)}$ 
of nonnegative integers
which is not $\Vec{0}$.
\item
A set $\Vec{\lambda}=(\lambda_e)_{e\in E_3(\Gamma)}$
of partitions
is called a {\em $\Gamma$-partition}.
$\Vec{\lambda}$ is {\em of degree} $\Vec{d}$ if
$(|\lambda_e|)_{e\in E_3(\Gamma)}=\Vec{d}$.
Note that
picking one $\Gamma$-partition is
the same as assigning a partition to every edge of $E_3(\Gamma)$.
\item
Given a $\Gamma$-partition $\Vec{\lambda}$,
we define 
$\Vec{\lambda}_v$ for a vertex $v\in V_3(\Gamma)$
as follows.
\begin{equation}\notag
\begin{split}
\\
&
\threein
\hspace*{2.2cm}
\twoin
\hspace*{2.2cm}
\onein
\hspace*{2.2cm}
\threeout
\\
&\Vec{\lambda}_v=(\lambda,\mu,\nu)\qquad
\Vec{\lambda}_v=(\lambda^t,\mu,\nu)\qquad
\Vec{\lambda}_v=(\lambda^t,\mu^t,\nu)\qquad
\Vec{\lambda}_v=(\lambda^t,\mu^t,\nu^t)
\end{split}
\end{equation}
It depends on
the directions of three incident edges and
their partitions.
If an incident edge is not in $E_3(\Gamma)$,
then we assume that the empty partition $\emptyset$ 
is assigned to it.
(Although  such edge is not directed,
it is irrelevant since $\emptyset^t=\emptyset$.)
\item
For a $\Gamma$-partition $\Vec{\lambda}$, we set
\begin{equation}\label{defy}
Y_{\Vec{\lambda}}(q)
=\prod_{e\in E_3(\Gamma)}
(-1)^{d_e(n_e+1)}q^{\frac{n_e\kappa(\lambda_e)}{2}}
\prod_{v\in V_3(\Gamma)}
C_{{\Vec{\lambda}}_v}(q).
\end{equation}
\end{itemize}

%
%
\begin{definition}%
The {\em partition function} of a GT graph $\Gamma$ is 
\begin{equation}\notag
{\mathcal Z}^{\Gamma}(q;\Vec{Q})=
1+\sum_{\Vec{d};\text{degree}}
{\mathcal Z}^{\Gamma}_{\Vec{d}}(q)\,
\Vec{Q}^{\Vec{d}}
,
\end{equation}
where $\Vec{Q}^{\Vec{d}}=\prod_{e\in E_3(\Gamma)}Q_e^{d_e}$
and
\begin{equation}\notag
{\mathcal Z}^{\Gamma}_{\Vec{d}}(q)
=\sum_{\begin{subarray}{c}
\Vec{\lambda}; \text{$\Gamma$-partition}\\
\text{of degree $\Vec{d}$}
\end{subarray}}
Y_{\Vec{\lambda}}(q).
\end{equation}
\end{definition}
\begin{definition}
The {\em free energy } of $\Gamma$ is defined as
\begin{equation}\notag
{\mathcal F}^{\Gamma}(q;\Vec{Q})=
\log {\mathcal Z}^{\Gamma}(q;\Vec{Q}).
\end{equation}
The coefficient of $\Vec{Q}^{\Vec{d}}$ is denoted by
${\mathcal F}_{\Vec{d}}^{\Gamma}(q)$.
\end{definition}

\subsection{Examples of Partition Function}\label{expart}
We calculate the partition function 
for the GT graphs in figure \ref{exgts}.
\subsubsection{}
%
First, we compute the partition function for the left GT graph.
Let us name trivalent vertices and the middle edge as follows.
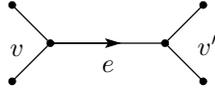
\begin{center}
\unitlength 0.1in
\begin{picture}( 10.1000,  4.0000)( 11.9000,-16.0000)
%
\special{pn 8}%
\special{pa 1200 1200}%
\special{pa 1400 1400}%
\special{fp}%
\special{pa 1400 1400}%
\special{pa 1200 1600}%
\special{fp}%
\special{pa 1400 1400}%
\special{pa 2000 1400}%
\special{fp}%
\special{pa 2000 1400}%
\special{pa 2200 1200}%
\special{fp}%
\special{pa 2000 1400}%
\special{pa 2200 1600}%
\special{fp}%
%
\special{pn 8}%
\special{pa 1400 1400}%
\special{pa 1760 1400}%
\special{fp}%
\special{sh 1}%
\special{pa 1760 1400}%
\special{pa 1694 1380}%
\special{pa 1708 1400}%
\special{pa 1694 1420}%
\special{pa 1760 1400}%
\special{fp}%
%
\special{pn 20}%
\special{sh 1}%
\special{ar 1400 1400 10 10 0  6.28318530717959E+0000}%
\special{sh 1}%
\special{ar 1200 1200 10 10 0  6.28318530717959E+0000}%
\special{sh 1}%
\special{ar 1200 1600 10 10 0  6.28318530717959E+0000}%
\special{sh 1}%
\special{ar 2200 1600 10 10 0  6.28318530717959E+0000}%
\special{sh 1}%
\special{ar 2200 1200 10 10 0  6.28318530717959E+0000}%
\special{sh 1}%
\special{ar 2000 1400 10 10 0  6.28318530717959E+0000}%
\put(16.7000,-15.5000){\makebox(0,0)[lb]{$e$}}%
\put(21.7000,-14.6000){\makebox(0,0)[lb]{$v'$}}%
\put(11.9000,-14.7000){\makebox(0,0)[lb]{$v$}}%
\end{picture}%
\end{center}
A $\Gamma$-partition consists of only one 
partition associated to the edge $e$:
$\Vec{\lambda}=(\lambda)$.
For this $\Gamma$-partition,
\begin{equation}\notag
\begin{split}
Y_{\Vec{\lambda}}(q)&=
(-1)^{(n+1)|\lambda|}q^{n\frac{\kappa(\lambda)}{2}}
C_{\lambda,\emptyset,\emptyset}(q)
C_{\lambda^t,\emptyset,\emptyset}(q)
\\
&=
(-1)^{(n+1)|\lambda|}
q^{n\frac{\kappa(\lambda)}{2}}
s_{\lambda}(q^{\rho})s_{\lambda^t}(q^{\rho})
\\
&=
(-1)^{(n+1)|\lambda|}
q^{(n-1)\frac{\kappa(\lambda)}{2}}
s_{\lambda}(q^{\rho})^2.
\end{split}
\end{equation}
In the last line, 
we have used the identity $s_{\lambda^t}(q^{\rho})
=q^{-\frac{\kappa(\lambda)}{2}}s_{\lambda}(q^{\rho})$
\cite{z2}.

Since a degree $\Vec{d}$ consists of only one component
$d$ associated to the edge $e$, 
we write $d$ instead of $\Vec{d}$.
We also write $Q_e$ as $Q$ for simplicity.
The partition function is
\begin{equation}\notag
\begin{split}
{\mathcal Z}^{\Gamma}(q;Q)&=
1+\sum_{d=1}^{\infty}{\mathcal Z}_{d}^{\Gamma}(q)\, Q^d,
\\
{\mathcal Z}_{d}^{\Gamma}(q)&=
(-1)^{(n+1)d}
\sum_{\lambda\in {\mathcal P}_d}
q^{(n-1)\frac{\kappa(\lambda)}{2}}
s_{\lambda}(q^{\rho})^2.
\end{split}
\end{equation}

This GT graph represents
the total space of 
${\mathcal O}(n-1)\oplus {\mathcal O}(-n-1)\to \mathbb{P}^1$
and
the free energy  ${\mathcal F}^{\Gamma}(q;Q)$
is nothing but the generating function of the GW invariants. 
%
\subsubsection{}
Next, we compute the partition function for the middle GT graph.
We introduce the {\em two-point function}
\begin{equation}\notag
W_{\mu,\nu}(q)=
(-1)^{|\mu|+|\nu|}q^{\frac{\kappa(\mu)+\kappa(\nu)}{2}}
\sum_{\eta\in{\mathcal P}}
s_{\mu/\eta}(q^{-\rho})
s_{\nu/\eta}(q^{-\rho})
\qquad (\mu,\nu\in{\mathcal P}),
\end{equation}
where $q^{-\rho}=(q^{i-\frac{1}{2}})_{i\geq 1}$.
It is a rational function in $q^{\frac{1}{2}}$
and satisfies $q^{\frac{\kappa(\mu)}{2}}W_{\mu,\nu}(q)=
C_{\mu^t,\emptyset,\nu}(q)$ (proposition 4.5, \cite{z2}).
Let us name the edge with the framing $\gamma_i+1$ as $e_i$
$(1\leq i\leq r)$
and the trivalent vertex incident on $e_i$ and $e_{i+1}$ as $v_i$.
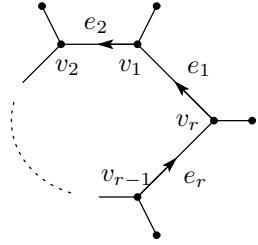
\begin{center}
\unitlength 0.1in
\begin{picture}( 12.7000, 12.3000)( 35.4000,-22.0000)
%
\special{pn 8}%
\special{pa 4000 1200}%
\special{pa 4200 1200}%
\special{fp}%
\special{pa 4200 1200}%
\special{pa 4600 1600}%
\special{fp}%
\special{pa 4600 1600}%
\special{pa 4200 2000}%
\special{fp}%
\special{pa 4200 2000}%
\special{pa 4000 2000}%
\special{fp}%
%
\special{pn 8}%
\special{pa 4010 1200}%
\special{pa 3800 1200}%
\special{fp}%
\special{pa 3800 1200}%
\special{pa 3600 1400}%
\special{fp}%
%
\special{pn 8}%
\special{ar 4000 1600 460 400  1.9028558 1.9307628}%
\special{ar 4000 1600 460 400  2.0144837 2.0423907}%
\special{ar 4000 1600 460 400  2.1261116 2.1540186}%
\special{ar 4000 1600 460 400  2.2377395 2.2656465}%
\special{ar 4000 1600 460 400  2.3493674 2.3772744}%
\special{ar 4000 1600 460 400  2.4609953 2.4889023}%
\special{ar 4000 1600 460 400  2.5726232 2.6005302}%
\special{ar 4000 1600 460 400  2.6842511 2.7121581}%
\special{ar 4000 1600 460 400  2.7958791 2.8237860}%
\special{ar 4000 1600 460 400  2.9075070 2.9354139}%
\special{ar 4000 1600 460 400  3.0191349 3.0470418}%
\special{ar 4000 1600 460 400  3.1307628 3.1586697}%
\special{ar 4000 1600 460 400  3.2423907 3.2702977}%
\special{ar 4000 1600 460 400  3.3540186 3.3819256}%
%
\special{pn 8}%
\special{pa 4600 1600}%
\special{pa 4810 1600}%
\special{fp}%
%
\special{pn 8}%
\special{pa 4200 2000}%
\special{pa 4300 2200}%
\special{fp}%
%
\special{pn 8}%
\special{pa 4300 1000}%
\special{pa 4200 1200}%
\special{fp}%
%
\special{pn 8}%
\special{pa 3800 1200}%
\special{pa 3700 1000}%
\special{fp}%
%
\special{pn 20}%
\special{sh 1}%
\special{ar 3800 1200 10 10 0  6.28318530717959E+0000}%
\special{sh 1}%
\special{ar 4200 1200 10 10 0  6.28318530717959E+0000}%
\special{sh 1}%
\special{ar 4600 1600 10 10 0  6.28318530717959E+0000}%
\special{sh 1}%
\special{ar 4200 2000 10 10 0  6.28318530717959E+0000}%
\special{sh 1}%
\special{ar 4300 2200 10 10 0  6.28318530717959E+0000}%
\special{sh 1}%
\special{ar 4800 1600 10 10 0  6.28318530717959E+0000}%
\special{sh 1}%
\special{ar 4300 1000 10 10 0  6.28318530717959E+0000}%
\special{sh 1}%
\special{ar 3700 1000 10 10 0  6.28318530717959E+0000}%
%
\special{pn 8}%
\special{pa 4200 2000}%
\special{pa 4400 1800}%
\special{fp}%
\special{sh 1}%
\special{pa 4400 1800}%
\special{pa 4340 1834}%
\special{pa 4362 1838}%
\special{pa 4368 1862}%
\special{pa 4400 1800}%
\special{fp}%
%
\special{pn 8}%
\special{pa 4600 1600}%
\special{pa 4400 1400}%
\special{fp}%
\special{sh 1}%
\special{pa 4400 1400}%
\special{pa 4434 1462}%
\special{pa 4438 1438}%
\special{pa 4462 1434}%
\special{pa 4400 1400}%
\special{fp}%
%
\special{pn 8}%
\special{pa 4200 1200}%
\special{pa 4000 1200}%
\special{fp}%
\special{sh 1}%
\special{pa 4000 1200}%
\special{pa 4068 1220}%
\special{pa 4054 1200}%
\special{pa 4068 1180}%
\special{pa 4000 1200}%
\special{fp}%
\put(44.4000,-19.5000){\makebox(0,0)[lb]{$e_r$}}%
\put(44.6000,-13.7000){\makebox(0,0)[lb]{$e_1$}}%
\put(39.2000,-11.4000){\makebox(0,0)[lb]{$e_2$}}%
\put(40.1000,-19.5000){\makebox(0,0)[lb]{$v_{r-1}$}}%
\put(44.0000,-16.5000){\makebox(0,0)[lb]{$v_r$}}%
\put(41.0000,-13.6000){\makebox(0,0)[lb]{$v_1$}}%
\put(37.7000,-13.6000){\makebox(0,0)[lb]{$v_2$}}%
\end{picture}%
\end{center}
%
%
Let $\lambda=(\lambda^1,\ldots,\lambda^r)$ 
be a $\Gamma$-partition where 
$\lambda^i$ is a partition assigned to edge $e_i$ $(1\leq i\leq r)$.
For $v_i$, $\Vec{\lambda}_{v_i}=(\lambda^{it},\emptyset,\lambda^{i+1})$.
Therefore
\begin{equation}\notag
\begin{split}
Y_{\Vec{\lambda}}(q)&=
\prod_{i=1}^r 
(-1)^{\gamma_i|\lambda^i|}
q^{(\gamma_i+1)\frac{\kappa(\lambda^i)}{2}}
C_{\lambda^{it},\emptyset,\lambda^{i+1}}(q)
\\
&=
\prod_{i=1}^r(-1)^{\gamma_i|\lambda^i|}q^{\frac{\gamma_i\kappa(\lambda^i)}{2}}
W_{\lambda^i,\lambda^{i+1}}(q).
\end{split}
\end{equation}
Here
$\lambda^{r+1}=\lambda^1$ is assumed.

We associate formal variables $Q_1,\ldots, Q_r$ to 
$e_1,\ldots, e_r$, respectively (In the previous notation,
$Q_i=Q_{e_i}$).
Then the partition function is
\begin{equation}\notag
{\mathcal Z}^{\Gamma}(q;Q_1,\ldots, Q_r)=
1+\sum_{\begin{subarray}{c}\Vec{d}=(d_1,\ldots,d_r);\\
\Vec{d}\neq \Vec{0},\\
d_i\geq 0\end{subarray}}
\prod_{i=1}^r
(-1)^{\gamma_i d_i}Q_i^{d_i}
\sum_{\begin{subarray}{c}(\lambda^1,\ldots,\lambda^r)\\
\lambda^i\in{\mathcal P_{d_i}}
\end{subarray}}
\prod_{i=1}^rq^{\frac{\gamma_i\kappa(\lambda^i)}{2}}
W_{\lambda^i,\lambda^{i+1}}(q).
\end{equation}

The GT graph represents a complete smooth toric surface $S$
if
$(\gamma_1,\ldots,\gamma_r)$
is equal to the set of self-intersection numbers
of the toric invariant curves in $S$.
In such a case,
the free energy ${\mathcal F}^{\Gamma}(q,\Vec{Q})$
is equal to the generating function of 
the GW invariants of the canonical bundle of $S$.
%
\subsubsection{}
Finally, we compute the partition function of 
the right GT graph.
Let us name trivalent vertices and edge
as follows.

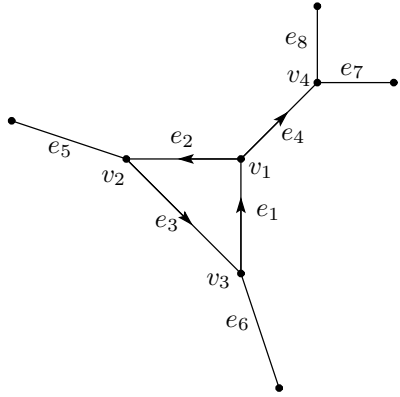
\begin{center}
\unitlength 0.1in
\begin{picture}( 20.0000, 20.0000)( 20.0000,-30.0000)
%
\special{pn 8}%
\special{pa 2600 1800}%
\special{pa 3200 1800}%
\special{pa 3200 2400}%
\special{pa 2600 1800}%
\special{fp}%
%
\special{pn 8}%
\special{pa 3200 2400}%
\special{pa 3400 3000}%
\special{fp}%
%
\special{pn 8}%
\special{pa 2000 1600}%
\special{pa 2600 1800}%
\special{fp}%
%
\special{pn 8}%
\special{pa 3600 1400}%
\special{pa 3200 1800}%
\special{fp}%
\special{pa 3600 1000}%
\special{pa 3600 1400}%
\special{fp}%
\special{pa 3600 1400}%
\special{pa 4000 1400}%
\special{fp}%
%
\special{pn 8}%
\special{pa 3200 2400}%
\special{pa 3200 2010}%
\special{fp}%
\special{sh 1}%
\special{pa 3200 2010}%
\special{pa 3180 2078}%
\special{pa 3200 2064}%
\special{pa 3220 2078}%
\special{pa 3200 2010}%
\special{fp}%
%
\special{pn 8}%
\special{pa 3200 1800}%
\special{pa 2880 1800}%
\special{fp}%
\special{sh 1}%
\special{pa 2880 1800}%
\special{pa 2948 1820}%
\special{pa 2934 1800}%
\special{pa 2948 1780}%
\special{pa 2880 1800}%
\special{fp}%
%
\special{pn 8}%
\special{pa 2600 1800}%
\special{pa 2940 2140}%
\special{fp}%
\special{sh 1}%
\special{pa 2940 2140}%
\special{pa 2908 2080}%
\special{pa 2902 2102}%
\special{pa 2880 2108}%
\special{pa 2940 2140}%
\special{fp}%
%
\special{pn 8}%
\special{pa 3200 1800}%
\special{pa 3440 1560}%
\special{fp}%
\special{sh 1}%
\special{pa 3440 1560}%
\special{pa 3380 1594}%
\special{pa 3402 1598}%
\special{pa 3408 1622}%
\special{pa 3440 1560}%
\special{fp}%
%
\special{pn 20}%
\special{sh 1}%
\special{ar 2000 1600 10 10 0  6.28318530717959E+0000}%
\special{sh 1}%
\special{ar 2600 1800 10 10 0  6.28318530717959E+0000}%
\special{sh 1}%
\special{ar 3200 2400 10 10 0  6.28318530717959E+0000}%
\special{sh 1}%
\special{ar 3400 3000 10 10 0  6.28318530717959E+0000}%
\special{sh 1}%
\special{ar 3200 1800 10 10 0  6.28318530717959E+0000}%
\special{sh 1}%
\special{ar 3600 1400 10 10 0  6.28318530717959E+0000}%
\special{sh 1}%
\special{ar 3600 1000 10 10 0  6.28318530717959E+0000}%
\special{sh 1}%
\special{ar 4000 1400 10 10 0  6.28318530717959E+0000}%
\put(32.8000,-21.2000){\makebox(0,0)[lb]{$e_1$}}%
\put(28.3000,-17.4000){\makebox(0,0)[lb]{$e_2$}}%
\put(27.5000,-21.8000){\makebox(0,0)[lb]{$e_3$}}%
\put(34.1000,-17.2000){\makebox(0,0)[lb]{$e_4$}}%
\put(21.9000,-17.9000){\makebox(0,0)[lb]{$e_5$}}%
\put(31.2000,-27.0000){\makebox(0,0)[lb]{$e_6$}}%
\put(37.2000,-13.8000){\makebox(0,0)[lb]{$e_7$}}%
\put(34.3000,-12.2000){\makebox(0,0)[lb]{$e_8$}}%
\put(32.4000,-19.0000){\makebox(0,0)[lb]{$v_1$}}%
\put(24.7000,-19.4000){\makebox(0,0)[lb]{$v_2$}}%
\put(30.2000,-24.8000){\makebox(0,0)[lb]{$v_3$}}%
\put(34.4000,-14.2000){\makebox(0,0)[lb]{$v_4$}}%
\end{picture}%

\end{center}
$E_3(\Gamma)=\{e_1,e_2,e_3,e_4\}$ and 
$V_3(\Gamma)=\{v_1,v_2,v_3,v_4\}$.
Let $\Vec{\lambda}=(\lambda^1,\ldots,\lambda^4)$ be
a $\Gamma$-partition where
$\lambda^i$ is a partition assigned to the edge $e_i$.
For each trivalent vertex,
\begin{equation}\label{ord}
\begin{split}
&\Vec{\lambda}_{v_1}=(\lambda^{1t},\lambda^4,\lambda^2),
\qquad
\Vec{\lambda}_{v_2}=(\lambda^{2t},\emptyset,\lambda^3),
\\
&
\Vec{\lambda}_{v_3}=(\lambda^{3t},\emptyset,\lambda^1),
\qquad
\Vec{\lambda}_{v_4}=(\emptyset,\lambda^{4t},\emptyset).
\end{split}
\end{equation}
Therefore
\begin{equation}\notag
Y_{\Vec{\lambda}}(q)=
(-1)^{\sum_{i=1}^4(b_i+1)|\lambda^i|}
q^{\sum_{i=1}^4 b_i\frac{\kappa(\lambda^i)}{2}}
C_{\lambda^{1t},\lambda^4,\lambda^{2t}}(q)
C_{\lambda^{2t},\emptyset,\lambda^{3t}}(q)
C_{\lambda^{3t},\emptyset,\lambda^{1t}}(q)
C_{\emptyset,\lambda^{4t},\emptyset}
\end{equation}
and
the partition function is
\begin{equation}\notag
\begin{split}
{\mathcal Z}^{\Gamma}(q;\Vec{Q})&=1+
\sum_{\begin{subarray}{c}
\Vec{d}=(d_1,d_2,d_3,d_4);\\
\Vec{d}\neq \Vec{0},\\
d_i\geq 0
\end{subarray}}
{\mathcal Z}_{\Vec{d}}^{\Gamma}(q)\,
Q_1^{d_1}\ldots Q_4^{d_4},
\\
{\mathcal Z}_{\Vec{d}}^{\Gamma}(q)&=
\sum_{\begin{subarray}{c}
\Vec{\lambda}=(\lambda^1,\lambda^2,\lambda^3,\lambda^4);\\
\lambda^i\in {\mathcal P}_{d_i}
\end{subarray}}
Y_{\Vec{\lambda}}(q).
\end{split}
\end{equation}

When $b_1=b_2=b_3=2$ and $b_4=0$,
the GT graph represents the flop of
the total space of the 
canonical bundle of the Hirzebruch surface $\mathbb{F}_1$
and the free energy is equal to the 
generating function of the GW invariants.
\section{Main Results}  \label{mainresults}
In this section, we state main results of this article.
Let us define
\begin{definition}\label{defG}
\begin{equation}\notag
G_{\Vec{ d}}^{\Gamma}(q)=
\sum_{k;k|d_0}\frac{\mu(k)}{k}
{\mathcal F}_{{\Vec{ d}}/k}^{\Gamma}(q^{k})\qquad
(d_0={\rm gcd}(\Vec{d})),
\end{equation}
where $\mu(k)$ is the M\"obius function.
\end{definition} %
%
We set $t=[1]^2$ and define ${\mathcal L}[t]$ by
\begin{equation}\notag
{\mathcal L}[t]=\Big\{\frac{f_2(t)}{f_1(t)}\Big|
f_1(t),f_2(t)\in \mathbb{Z}[t],f_1(t):\text{ monic}
\Big\}.
\end{equation}
${\mathcal L}[t]$ is a subring of the ring of rational functions
$\mathbb{Q}(t)$ \cite{peng}.

The main results of the paper are
\begin{prop}\label{integrality}
\begin{equation}\notag
G_{\Vec{d}}^{\Gamma}(q)\in{\mathcal L}[t].
\end{equation}
\end{prop}
\begin{prop}\label{polynomiality}
\begin{equation}\notag
t \cdot G_{\Vec{d}}^{\Gamma}(q)\in\mathbb{Q}[t].
\end{equation}
\end{prop}
We  will prove propositions \ref{integrality}
and \ref{polynomiality} in sections \ref{PI} and \ref{PV},
respectively.

Propositions \ref{integrality} and
\ref{polynomiality}
imply that the 
numerator of $t\cdot G^{\Gamma}_{\Vec{d}}(q)$
is divisible by the denominator.
Since the denominator is monic,
the quotient is a polynomial in $t$ with integer
coefficients.
Thus
\begin{cor} \label{mainthm}
$t \cdot G_{\Vec{d}}^{\Gamma}(q)\in\mathbb{Z}[t]$.
\end{cor}%
%

What does this corollary mean ?
By the formula of the M\"obius function
\begin{equation}\label{keyformula}
\sum_{k':k'|k}\mu(k')
=\begin{cases}
1&(k=1)\\
0&(k>1,k\in\mathbb{N}),
\end{cases}
\end{equation}
the free energy in degree $\Vec{d}$ is written as
\begin{equation}\notag
{\mathcal F}_{\Vec{d}}^{\Gamma}(q)
=
\sum_{k;k|d_0}\frac{1}{k}
{G}_{{\Vec{ d}}/k}^{\Gamma}(q^{k}).
\end{equation}
In fact, 
definition  \ref{defG} was 
obtained by inverting this relation \cite{bp}.
Let us write the corollary as follows.
\begin{equation}\notag
G_{\Vec{d}}^{\Gamma}(q)=\sum_{g\geq 0}
n_{\Vec{d}}^g(\Gamma)(-t)^{g-1}
\end{equation}
where $\{n_{\Vec{d}}^g(\Gamma)\}_{g\geq 0}$
is 
a sequence of integers 
only finite number of which is nonzero.
Note that
proposition \ref{integrality}  implies the 
integrality of $n_{\Vec{d}}^g(\Gamma)$.
Proposition
\ref{polynomiality} implies the
vanishing of $n_{\Vec{d}}^g$ at large $g$
(and also at $g<0$).
We find that the free energy is written in
terms of these integers as
\begin{equation}\label{fn}
{\mathcal F}_{\Vec{d}}^{\Gamma}(q)
=\sum_{g\geq 0}\sum_{k;k|d_0}
n_{\Vec{d}/k}^g(\Gamma)
\frac{(-t_k)^{g-1}}{k},
\end{equation}
where $t_k=[k]^2$.
%
%
%

%
Before moving to the proof of the propositions,
we explain the geometric meaning of these results. 

\section{Toric Calabi--Yau 
Threefold and Gopakumar--Vafa conjecture}\label{TCY}
Given a toric Calabi--Yau threefold (TCY threefold) $X$, 
a planar graph is 
determined canonically from the fan of $X$.
It is called the {\em toric graph} of $X$ and it is
a GT graph or the graph union of GT graphs.
In this section,
we first describe how to draw the toric graph.
Then we explain the relation between the free energy of 
the toric graph 
and the generating function of the GW invariants of $X$.
Finally,
we see that (\ref{fn})
implies the integrality and
the vanishing for high genera of the 
Gopakumar--Vafa invariants.
%
\subsection{TCY threefold}
A {\em Calabi--Yau toric threefold} is 
a three-dimensional, smooth toric variety $X$
of finite type,
whose canonical bundle $K_X$ is a trivial line bundle.
The last condition is called the 
{\em Calabi--Yau condition}.
For simplicity, we impose one more condition,
which implies 
that the fundamental group $\pi_1(X)$ is trivial 
and that $H^2(X)\cong {\rm Pic}(X)$.

A toric variety $X$ is constructed from a fan $\Sigma$,
which is a collection of cones.
The fan of $X$ is unique up to $SL(3,\mathbb{Z})$
since the action of $SL(3,\mathbb{Z})$
on a fan is offset by the change of the coordinate functions.

The conditions on $X$ is rephrased in terms of those on
the fan $\Sigma$ as follows.
\begin{description}
\item[Finite type]
$X$ is of finite type if its fan $\Sigma$ is a finite set.
\item[Smoothness]
$X$ is smooth if and only if
the minimal set of generators of every cone 
forms a part of a $\mathbb{Z}$-basis 
of $\mathbb{R}^3$.
(Here the generators of a cone mean the
shortest integral vectors that generate the cone.)
%
\item[Calabi--Yau]
The canonical bundle of $X$ is trivial if and only if
there exists a vector $u\in(\mathbb{R}^3)^*$ satisfying
\begin{equation}\notag
\langle \omega_i,u\rangle =1
\end{equation}
for all generators $\omega_i$ of the fan.
Using the action of $SL(3,\mathbb{Z})$,
we take 
\begin{equation}\notag
u=(0,0,1).
\end{equation}
Therefore every generators of a fan of 
a toric Calabi--Yau threefolds
is of the form $(*,*,1)$.
Note that such  fan can not be complete.
Equivalently, the toric variety $X$ is noncompact.
\item[Other assumption]
We assume that 
there exists at least one 3-cone and that
every 1 or 2-cone of the fan $\Sigma$ is a face of some 3-cone.
This
implies that
\begin{equation}\notag
\pi_1(X)=\{id\},\qquad H^2(X)\cong {\rm Pic}(X).
\end{equation}
See \cite{fulton} for a proof.
\end{description}

\subsection{Toric Graph}
Since all the generators are of the form $(*,*,1)$,
it is sufficient 
to see the section $\bar{\Sigma}$ of
the fan $\Sigma$ at the height 1.
%
We will write the section of a cone $\sigma$ as $\bar{\sigma}$.

From $\bar{\Sigma}$,
we draw a labeled graph as follows.
\begin{enumerate}
\item
Draw a vertex $v_{\sigma}$ inside every 2-simplex $\bar{\sigma}$.
%
\item
Draw an edge $e_{\tau}$ transversally to every 1-simplex 
$\bar{\tau}$ as follows.
\begin{enumerate}
\item \label{Ethree}
If $\bar{\tau}$ is the boundary of 
two 2-simplices $\bar{\sigma},\bar{\sigma}'$,
let $e_{\tau}$ join $v_{\sigma}$ and $v_{\sigma'}$.
%
%
\item
If $\bar{\tau}$ is the boundary of 
only one 2-simplex $\bar{\sigma}$,
let $e_{\tau}$ be incident to $v_{\sigma}$;
add one vertex $v_{\tau}$ to other endpoint.
\end{enumerate}
%
%
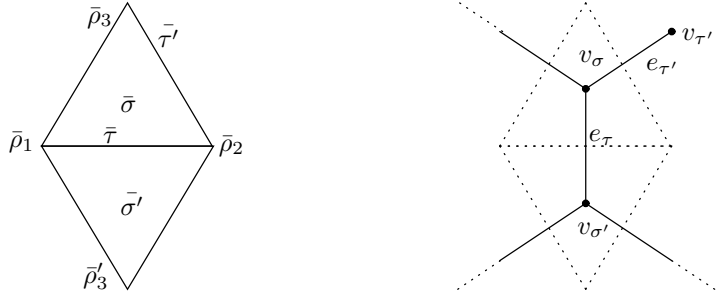
\begin{center}
\unitlength 0.1in
\begin{picture}( 36.9600, 15.4000)( 15.7000,-25.0000)
%
\special{pn 8}%
\special{pa 1742 1750}%
\special{pa 2636 1750}%
\special{pa 2192 1000}%
\special{pa 1742 1750}%
\special{fp}%
%
\special{pn 8}%
\special{pa 1742 1750}%
\special{pa 2192 2500}%
\special{pa 2642 1750}%
\special{pa 1742 1750}%
\special{fp}%
%
\special{pn 8}%
\special{pa 4140 1750}%
\special{pa 5034 1750}%
\special{pa 4590 1000}%
\special{pa 4140 1750}%
\special{dt 0.045}%
%
\special{pn 8}%
\special{pa 4140 1750}%
\special{pa 4590 2500}%
\special{pa 5040 1750}%
\special{pa 4140 1750}%
\special{dt 0.045}%
\put(21.5500,-15.7000){\makebox(0,0)[lb]{$\bar{\sigma}$}}%
\put(21.5500,-21.1000){\makebox(0,0)[lb]{$\bar{\sigma'}$}}%
\put(20.6500,-17.3500){\makebox(0,0)[lb]{$\bar{\tau}$}}%
\put(23.5000,-12.3200){\makebox(0,0)[lb]{$\bar{\tau'}$}}%
%
\special{pn 8}%
\special{pa 4590 1450}%
\special{pa 4590 2050}%
\special{fp}%
%
\special{pn 8}%
\special{pa 4590 2050}%
\special{pa 5040 2350}%
\special{fp}%
%
\special{pn 8}%
\special{pa 4140 1150}%
\special{pa 4590 1450}%
\special{fp}%
%
\special{pn 8}%
\special{pa 4590 1450}%
\special{pa 5040 1150}%
\special{fp}%
%
\special{pn 8}%
\special{pa 4140 2350}%
\special{pa 4590 2050}%
\special{fp}%
%
\special{pn 8}%
\special{pa 5040 2350}%
\special{pa 5266 2500}%
\special{dt 0.045}%
%
\special{pn 8}%
\special{pa 3930 1000}%
\special{pa 4156 1150}%
\special{dt 0.045}%
%
\special{pn 8}%
\special{pa 3920 2500}%
\special{pa 4144 2350}%
\special{dt 0.045}%
%
\special{pn 20}%
\special{sh 1}%
\special{ar 5040 1150 10 10 0  6.28318530717959E+0000}%
\special{sh 1}%
\special{ar 4590 1450 10 10 0  6.28318530717959E+0000}%
\special{sh 1}%
\special{ar 4590 2050 10 10 0  6.28318530717959E+0000}%
\put(50.9300,-12.1800){\makebox(0,0)[lb]{$v_{\tau'}$}}%
\put(45.5300,-13.3000){\makebox(0,0)[lb]{$v_{\sigma}$}}%
\put(45.4600,-22.3000){\makebox(0,0)[lb]{$v_{\sigma'}$}}%
\put(46.0600,-17.4200){\makebox(0,0)[lb]{$e_{\tau}$}}%
\put(15.7000,-17.8000){\makebox(0,0)[lb]{$\bar{\rho}_1$}}%
\put(26.7000,-18.0000){\makebox(0,0)[lb]{$\bar{\rho}_2$}}%
\put(49.0600,-13.7500){\makebox(0,0)[lb]{$e_{\tau'}$}}%
\put(19.8000,-11.3000){\makebox(0,0)[lb]{$\bar{\rho}_3$}}%
\put(19.6000,-25.0000){\makebox(0,0)[lb]{$\bar{\rho}_3'$}}%
\end{picture}%
\end{center}

\item
To every flag $f$ whose edge is of type \ref{Ethree},
we assign an integer label $n_f$ as follows.
%
%
%
For $(v,\sigma)$ and $(v,\sigma')$ in the above figure, 
the labels are
\begin{equation}\notag
\frac{a_1-a_2}{2}\text{ for }(v_{\sigma},e_{\tau}),\qquad
\frac{-a_1+a_2}{2}\text{ for }(v_{\sigma'},e_{\tau}).
\end{equation}
Here $a_1,a_2$ are integers defined by
\begin{equation}\notag
\omega_3'=-a_1\omega_1-a_2\omega_2-\omega_3
\end{equation}
where $\omega_1,\omega_2,\omega_3$ and $\omega_3'$ are
generators of the 1-cones 
$\rho_1,\rho_2,\rho_3$ and $\rho_3'$, respectively.
Since $a_1+a_2=-2$ by the Calabi--Yau condition,
these are integers.
The label is called the {\em framing} of the flag.
For reference, we computed the framings for flags
in figure \ref{exframings}.
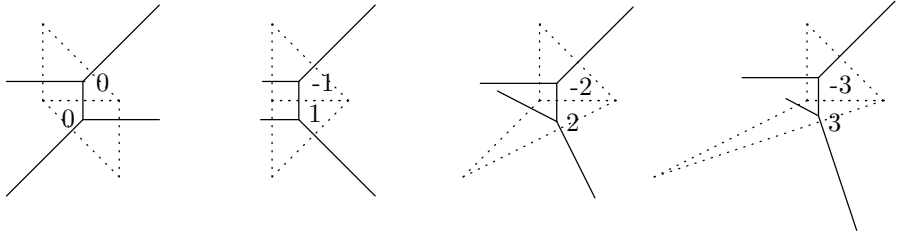
\begin{figure}[th]
\begin{center}
\unitlength 0.1in
\begin{picture}( 46.5000, 12.0000)(  8.1000,-20.8000)
%
\special{pn 8}%
\special{pa 1000 1400}%
\special{pa 1400 1400}%
\special{pa 1000 1000}%
\special{pa 1000 1000}%
\special{pa 1000 1400}%
\special{dt 0.045}%
%
\special{pn 8}%
\special{pa 2200 1400}%
\special{pa 2600 1400}%
\special{pa 2200 1000}%
\special{pa 2200 1000}%
\special{pa 2200 1400}%
\special{dt 0.045}%
%
\special{pn 8}%
\special{pa 3600 1400}%
\special{pa 4000 1400}%
\special{pa 3600 1000}%
\special{pa 3600 1000}%
\special{pa 3600 1400}%
\special{dt 0.045}%
%
\special{pn 8}%
\special{pa 1000 1400}%
\special{pa 1400 1800}%
\special{pa 1400 1400}%
\special{pa 1000 1400}%
\special{dt 0.045}%
%
\special{pn 8}%
\special{pa 2200 1400}%
\special{pa 2200 1800}%
\special{pa 2600 1400}%
\special{pa 2200 1400}%
\special{dt 0.045}%
%
\special{pn 8}%
\special{pa 5000 1400}%
\special{pa 5400 1400}%
\special{pa 5000 1000}%
\special{pa 5000 1000}%
\special{pa 5000 1400}%
\special{dt 0.045}%
%
\special{pn 8}%
\special{pa 3600 1400}%
\special{pa 3200 1800}%
\special{pa 4000 1400}%
\special{pa 3600 1400}%
\special{dt 0.045}%
%
\special{pn 8}%
\special{pa 5000 1400}%
\special{pa 4200 1800}%
\special{pa 5400 1400}%
\special{pa 5400 1400}%
\special{pa 5000 1400}%
\special{dt 0.045}%
%
\special{pn 8}%
\special{pa 1210 1300}%
\special{pa 1210 1500}%
\special{fp}%
\special{pa 1210 1500}%
\special{pa 810 1900}%
\special{fp}%
\special{pa 1210 1300}%
\special{pa 1610 900}%
\special{fp}%
\special{pa 1210 1500}%
\special{pa 1610 1500}%
\special{fp}%
\special{pa 1220 1300}%
\special{pa 810 1300}%
\special{fp}%
%
\special{pn 8}%
\special{pa 2150 1300}%
\special{pa 2340 1300}%
\special{fp}%
\special{pa 2340 1300}%
\special{pa 2340 1500}%
\special{fp}%
\special{pa 2340 1500}%
\special{pa 2140 1500}%
\special{fp}%
\special{pa 2340 1500}%
\special{pa 2740 1900}%
\special{fp}%
\special{pa 2340 1300}%
\special{pa 2740 900}%
\special{fp}%
%
\special{pn 8}%
\special{pa 3290 1310}%
\special{pa 3690 1310}%
\special{fp}%
\special{pa 3690 1310}%
\special{pa 4090 910}%
\special{fp}%
%
\special{pn 8}%
\special{pa 4660 1280}%
\special{pa 5060 1280}%
\special{fp}%
\special{pa 5060 1280}%
\special{pa 5460 880}%
\special{fp}%
%
\special{pn 8}%
\special{pa 3690 1310}%
\special{pa 3690 1510}%
\special{fp}%
\special{pa 3690 1510}%
\special{pa 3380 1350}%
\special{fp}%
\special{pa 3690 1510}%
\special{pa 3890 1910}%
\special{fp}%
%
\special{pn 8}%
\special{pa 5060 1280}%
\special{pa 5060 1480}%
\special{fp}%
%
\special{pn 8}%
\special{pa 5060 1480}%
\special{pa 4890 1390}%
\special{fp}%
\special{pa 5060 1480}%
\special{pa 5260 2080}%
\special{fp}%
\put(12.8000,-13.5000){\makebox(0,0)[lb]{0}}%
\put(11.0000,-15.4000){\makebox(0,0)[lb]{0}}%
\put(24.1000,-13.5000){\makebox(0,0)[lb]{-1}}%
\put(23.9000,-15.1000){\makebox(0,0)[lb]{1}}%
\put(37.6000,-13.7000){\makebox(0,0)[lb]{-2}}%
\put(37.4000,-15.6000){\makebox(0,0)[lb]{2}}%
\put(51.2000,-13.6000){\makebox(0,0)[lb]{-3}}%
\put(51.1000,-15.6000){\makebox(0,0)[lb]{3}}%
\end{picture}%
\end{center}
\caption{Examples of 
framings.
}\label{exframings}
\end{figure}
\end{enumerate}
The resulting graph is the {\em toric graph} of 
the TCY threefold $X$.
Note that the toric graph 
is unique
although the fan  is unique only up to the action of 
$SL(3,\mathbb{Z})$.

\begin{figure}[th]
\begin{center}
\unitlength 0.1in
\begin{picture}( 50.3800, 10.3600)(  9.6000,-19.5600)
%
\special{pn 8}%
\special{pa 1548 1536}%
\special{pa 1688 1256}%
\special{pa 1968 1536}%
\special{pa 1548 1536}%
\special{fp}%
%
\special{pn 8}%
\special{pa 1128 1816}%
\special{pa 1968 1536}%
\special{pa 1548 1536}%
\special{pa 1548 1536}%
\special{pa 1128 1816}%
\special{fp}%
%
\special{pn 8}%
\special{pa 1548 1536}%
\special{pa 1688 1256}%
\special{pa 1128 1816}%
\special{pa 1548 1536}%
\special{fp}%
%
\special{pn 8}%
\special{pa 1548 1956}%
\special{pa 2108 1396}%
\special{fp}%
%
\special{pn 8}%
\special{pa 1548 1956}%
\special{pa 960 1760}%
\special{fp}%
%
\special{pn 8}%
\special{pa 1548 1956}%
\special{pa 1548 1536}%
\special{dt 0.045}%
%
\special{pn 8}%
\special{pa 1548 1536}%
\special{pa 1548 1116}%
\special{fp}%
%
\special{pn 8}%
\special{pa 1688 1256}%
\special{pa 1548 1956}%
\special{dt 0.045}%
%
\special{pn 8}%
\special{pa 1758 920}%
\special{pa 1688 1270}%
\special{fp}%
%
\special{pn 8}%
\special{pa 5268 1536}%
\special{pa 5408 1256}%
\special{pa 5688 1536}%
\special{pa 5268 1536}%
\special{fp}%
%
\special{pn 8}%
\special{pa 5688 1536}%
\special{pa 5828 1256}%
\special{pa 5408 1256}%
\special{pa 5408 1256}%
\special{pa 5688 1536}%
\special{fp}%
%
\special{pn 8}%
\special{pa 4848 1816}%
\special{pa 5688 1536}%
\special{pa 5268 1536}%
\special{pa 5268 1536}%
\special{pa 4848 1816}%
\special{fp}%
%
\special{pn 8}%
\special{pa 5268 1536}%
\special{pa 5408 1256}%
\special{pa 4848 1816}%
\special{pa 5268 1536}%
\special{fp}%
%
\special{pn 8}%
\special{pa 5268 1956}%
\special{pa 5828 1396}%
\special{fp}%
%
\special{pn 8}%
\special{pa 5268 1956}%
\special{pa 5268 1536}%
\special{dt 0.045}%
%
\special{pn 8}%
\special{pa 5268 1536}%
\special{pa 5268 1116}%
\special{fp}%
%
\special{pn 8}%
\special{pa 5408 1256}%
\special{pa 5268 1956}%
\special{dt 0.045}%
%
\special{pn 8}%
\special{pa 5828 1256}%
\special{pa 5268 1956}%
\special{dt 0.045}%
%
\special{pn 8}%
\special{pa 3408 1536}%
\special{pa 3548 1256}%
\special{pa 3828 1536}%
\special{pa 3408 1536}%
\special{fp}%
%
\special{pn 8}%
\special{pa 2988 1816}%
\special{pa 3828 1536}%
\special{pa 3408 1536}%
\special{pa 3408 1536}%
\special{pa 2988 1816}%
\special{fp}%
%
\special{pn 8}%
\special{pa 3408 1956}%
\special{pa 3968 1396}%
\special{fp}%
%
\special{pn 8}%
\special{pa 3408 1956}%
\special{pa 3408 1536}%
\special{dt 0.045}%
%
\special{pn 8}%
\special{pa 3408 1536}%
\special{pa 3408 1116}%
\special{fp}%
%
\special{pn 8}%
\special{pa 3548 1256}%
\special{pa 3408 1956}%
\special{dt 0.045}%
%
\special{pn 8}%
\special{pa 3618 920}%
\special{pa 3548 1270}%
\special{fp}%
%
\special{pn 8}%
\special{pa 5478 920}%
\special{pa 5408 1270}%
\special{fp}%
%
\special{pn 8}%
\special{pa 3408 1956}%
\special{pa 2820 1760}%
\special{fp}%
%
\special{pn 8}%
\special{pa 5268 1956}%
\special{pa 4680 1760}%
\special{fp}%
%
\special{pn 8}%
\special{pa 5830 1260}%
\special{pa 5998 1050}%
\special{fp}%
\put(10.7200,-19.3500){\makebox(0,0)[lb]{$w_3$}}%
\put(19.6100,-16.5500){\makebox(0,0)[lb]{$w_1$}}%
\put(36.1000,-12.8000){\makebox(0,0)[lb]{$w_2$}}%
\put(13.7000,-15.6000){\makebox(0,0)[lb]{$w_0$}}%
\put(59.1000,-13.6000){\makebox(0,0)[lb]{$w_4$}}%
\put(17.5000,-12.9000){\makebox(0,0)[lb]{$w_2$}}%
\put(38.7000,-16.7000){\makebox(0,0)[lb]{$w_1$}}%
\put(57.3000,-16.5000){\makebox(0,0)[lb]{$w_1$}}%
\put(29.6000,-19.6000){\makebox(0,0)[lb]{$w_3$}}%
\put(48.2000,-19.8000){\makebox(0,0)[lb]{$w_3$}}%
\put(54.6000,-12.5000){\makebox(0,0)[lb]{$w_2$}}%
\put(32.4000,-15.6000){\makebox(0,0)[lb]{$w_0$}}%
\put(50.9000,-15.7000){\makebox(0,0)[lb]{$w_0$}}%
\end{picture}%

\unitlength 0.1in
\begin{picture}( 54.2000, 18.2000)( 10.0000,-28.0000)
%
\special{pn 8}%
\special{pa 1540 1860}%
\special{pa 1940 1860}%
\special{pa 1940 2260}%
\special{pa 1540 1860}%
\special{fp}%
%
\special{pn 8}%
\special{pa 1940 2260}%
\special{pa 2140 2660}%
\special{fp}%
%
\special{pn 8}%
\special{pa 1150 1660}%
\special{pa 1550 1860}%
\special{fp}%
%
\special{pn 8}%
\special{pa 1800 2000}%
\special{pa 2600 2000}%
\special{pa 1800 1200}%
\special{pa 1800 2000}%
\special{dt 0.045}%
%
\special{pn 8}%
\special{pa 1800 2000}%
\special{pa 1000 2800}%
\special{pa 2600 2000}%
\special{pa 1800 2000}%
\special{dt 0.045}%
%
\special{pn 8}%
\special{pa 1800 2000}%
\special{pa 1800 1200}%
\special{pa 1000 2800}%
\special{pa 1000 2800}%
\special{pa 1800 2000}%
\special{dt 0.045}%
%
\special{pn 8}%
\special{pa 1940 1860}%
\special{pa 2440 1360}%
\special{fp}%
%
\special{pn 8}%
\special{pa 3600 2000}%
\special{pa 4400 2000}%
\special{pa 3600 1200}%
\special{pa 3600 2000}%
\special{dt 0.045}%
%
\special{pn 8}%
\special{pa 3600 2000}%
\special{pa 2800 2800}%
\special{pa 4400 2000}%
\special{pa 3600 2000}%
\special{dt 0.045}%
%
\special{pn 8}%
\special{pa 5400 2000}%
\special{pa 6200 2000}%
\special{pa 5400 1200}%
\special{pa 5400 2000}%
\special{dt 0.045}%
%
\special{pn 8}%
\special{pa 5400 2000}%
\special{pa 4600 2800}%
\special{pa 6200 2000}%
\special{pa 5400 2000}%
\special{dt 0.045}%
%
\special{pn 8}%
\special{pa 5400 2000}%
\special{pa 5400 1200}%
\special{pa 4600 2800}%
\special{pa 4600 2800}%
\special{pa 5400 2000}%
\special{dt 0.045}%
%
\special{pn 8}%
\special{pa 6200 2000}%
\special{pa 6200 1200}%
\special{pa 5400 1200}%
\special{pa 6200 2000}%
\special{dt 0.045}%
%
\special{pn 8}%
\special{pa 3350 1860}%
\special{pa 3750 1860}%
\special{fp}%
\special{pa 3750 1860}%
\special{pa 3750 2260}%
\special{fp}%
%
\special{pn 8}%
\special{pa 3470 1980}%
\special{pa 3750 2260}%
\special{fp}%
%
\special{pn 8}%
\special{pa 3750 1860}%
\special{pa 4250 1360}%
\special{fp}%
%
\special{pn 8}%
\special{pa 3750 2260}%
\special{pa 3950 2660}%
\special{fp}%
%
\special{pn 8}%
\special{pa 5120 1880}%
\special{pa 5520 1880}%
\special{pa 5520 2280}%
\special{pa 5120 1880}%
\special{fp}%
%
\special{pn 8}%
\special{pa 5520 2280}%
\special{pa 5720 2680}%
\special{fp}%
%
\special{pn 8}%
\special{pa 4730 1680}%
\special{pa 5130 1880}%
\special{fp}%
%
\special{pn 8}%
\special{pa 5520 1880}%
\special{pa 6020 1380}%
\special{fp}%
%
\special{pn 8}%
\special{pa 6020 1380}%
\special{pa 6420 1380}%
\special{fp}%
\special{pa 6020 980}%
\special{pa 6020 1380}%
\special{fp}%
%
\special{pn 20}%
\special{sh 1}%
\special{ar 1940 2260 10 10 0  6.28318530717959E+0000}%
\special{sh 1}%
\special{ar 1940 1860 10 10 0  6.28318530717959E+0000}%
\special{sh 1}%
\special{ar 1550 1860 10 10 0  6.28318530717959E+0000}%
%
\special{pn 20}%
\special{sh 1}%
\special{ar 2430 1360 10 10 0  6.28318530717959E+0000}%
\special{sh 1}%
\special{ar 2140 2640 10 10 0  6.28318530717959E+0000}%
\special{sh 1}%
\special{ar 1140 1660 10 10 0  6.28318530717959E+0000}%
%
\special{pn 20}%
\special{sh 1}%
\special{ar 3490 1980 10 10 0  6.28318530717959E+0000}%
\special{sh 1}%
\special{ar 3350 1860 10 10 0  6.28318530717959E+0000}%
\special{sh 1}%
\special{ar 3750 1860 10 10 0  6.28318530717959E+0000}%
\special{sh 1}%
\special{ar 3750 2250 10 10 0  6.28318530717959E+0000}%
\special{sh 1}%
\special{ar 3950 2650 10 10 0  6.28318530717959E+0000}%
\special{sh 1}%
\special{ar 4240 1370 10 10 0  6.28318530717959E+0000}%
%
\special{pn 20}%
\special{sh 1}%
\special{ar 5120 1880 10 10 0  6.28318530717959E+0000}%
\special{sh 1}%
\special{ar 5530 1890 10 10 0  6.28318530717959E+0000}%
\special{sh 1}%
\special{ar 5520 2280 10 10 0  6.28318530717959E+0000}%
\special{sh 1}%
\special{ar 5730 2680 10 10 0  6.28318530717959E+0000}%
\special{sh 1}%
\special{ar 4730 1670 10 10 0  6.28318530717959E+0000}%
\special{sh 1}%
\special{ar 6020 1380 10 10 0  6.28318530717959E+0000}%
\special{sh 1}%
\special{ar 6400 1380 10 10 0  6.28318530717959E+0000}%
\special{sh 1}%
\special{ar 6030 990 10 10 0  6.28318530717959E+0000}%
\put(20.1000,-19.7000){\makebox(0,0)[lb]{-2}}%
\put(19.9000,-22.9000){\makebox(0,0)[lb]{2}}%
\put(15.5000,-18.3000){\makebox(0,0)[lb]{-2}}%
\put(18.6000,-18.3000){\makebox(0,0)[lb]{2}}%
\put(14.6000,-20.4000){\makebox(0,0)[lb]{2}}%
\put(17.4000,-23.4000){\makebox(0,0)[lb]{-2}}%
\put(52.9000,-23.4000){\makebox(0,0)[lb]{-2}}%
\put(37.9000,-22.9000){\makebox(0,0)[lb]{2}}%
\put(55.9000,-22.9000){\makebox(0,0)[lb]{2}}%
\put(38.1000,-19.7000){\makebox(0,0)[lb]{-2}}%
\put(56.1000,-19.7000){\makebox(0,0)[lb]{-2}}%
\put(50.6000,-20.4000){\makebox(0,0)[lb]{2}}%
\put(51.5000,-18.3000){\makebox(0,0)[lb]{-2}}%
\put(54.1000,-18.4000){\makebox(0,0)[lb]{2}}%
\put(58.9000,-14.5000){\makebox(0,0)[lb]{0}}%
\put(55.8000,-17.6000){\makebox(0,0)[lb]{0}}%
\end{picture}%
\end{center}
\caption{Examples of fans (upstairs)
and toric graphs (downstairs).
%
%
The left is 
the canonical bundle of $\mathbb{P}^2$,
the middle is 
the total space of the vector bundle 
${\mathcal O}(1)\oplus
{\mathcal O}(-3)\to \mathbb{P}^1$ and
the left is 
the flop of the canonical bundle of 
the Hirzebruch surface $\mathbb{F}_1$.
$w_i$ ($0\leq i\leq 4$) are the generators:
$w_0=(0,0,1),w_1=(1,0,1),w_2=(0,1,1),w_3=(-1,-1,1)$
and $w_4=(1,1,1)$.
}
\label{exfan}
\end{figure}
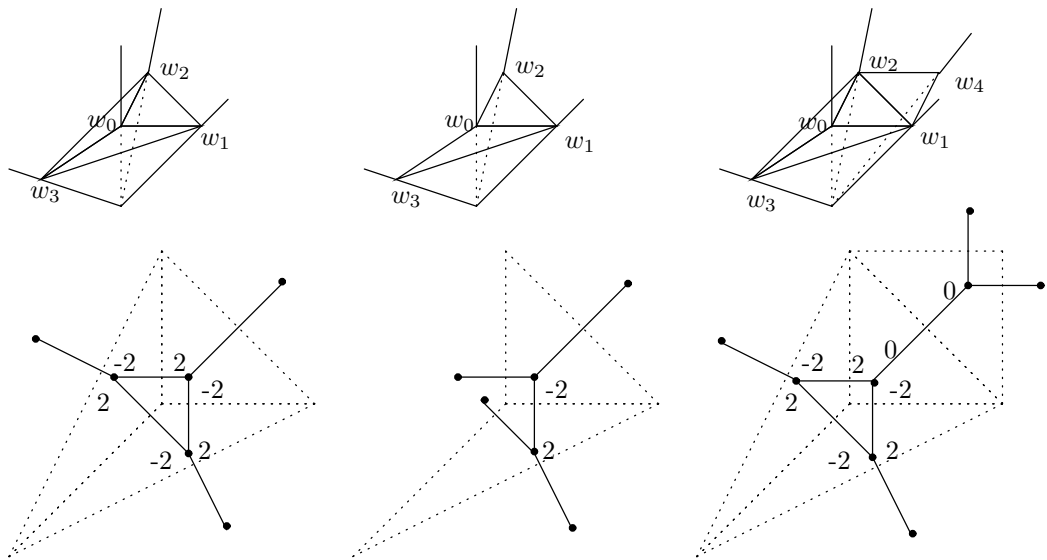

\begin{figure}[ht]
\begin{center}
\unitlength 0.1in
\begin{picture}(  8.0000,  8.0000)( 12.0000,-22.0000)
%
\special{pn 8}%
\special{pa 1600 1400}%
\special{pa 1600 1800}%
\special{pa 2000 1800}%
\special{pa 1600 1400}%
\special{dt 0.045}%
%
\special{pn 8}%
\special{pa 1200 1800}%
\special{pa 1600 1800}%
\special{pa 1600 2200}%
\special{pa 1200 1800}%
\special{dt 0.045}%
%
\special{pn 8}%
\special{pa 1750 1660}%
\special{pa 1750 1860}%
\special{fp}%
\special{pa 1750 1660}%
\special{pa 1550 1660}%
\special{fp}%
\special{pa 1750 1660}%
\special{pa 1950 1460}%
\special{fp}%
%
\special{pn 20}%
\special{sh 1}%
\special{ar 1750 1660 10 10 0  6.28318530717959E+0000}%
\special{sh 1}%
\special{ar 1750 1860 10 10 0  6.28318530717959E+0000}%
\special{sh 1}%
\special{ar 1560 1660 10 10 0  6.28318530717959E+0000}%
\special{sh 1}%
\special{ar 1950 1460 10 10 0  6.28318530717959E+0000}%
%
\special{pn 8}%
\special{pa 1450 1950}%
\special{pa 1650 1950}%
\special{fp}%
\special{pa 1450 1950}%
\special{pa 1450 1750}%
\special{fp}%
\special{pa 1450 1950}%
\special{pa 1250 2150}%
\special{fp}%
%
\special{pn 20}%
\special{sh 1}%
\special{ar 1450 1950 10 10 0  6.28318530717959E+0000}%
\special{sh 1}%
\special{ar 1650 1940 10 10 0  6.28318530717959E+0000}%
\special{sh 1}%
\special{ar 1450 1750 10 10 0  6.28318530717959E+0000}%
\special{sh 1}%
\special{ar 1250 2150 10 10 0  6.28318530717959E+0000}%
\end{picture}%
\end{center}
\caption{An example of 
the toric graph with more than one connected components.
This corresponds to
the canonical bundle of 
the noncomplete toric surface 
${\mathbb{P}^1}\times \mathbb{P}^1\setminus
\{(0,0),(\infty,\infty)\}$.}
\label{noncomplete}
\end{figure}
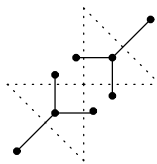

Examples of the toric graphs are shown in figures \ref{exfan}
and \ref{noncomplete}.
See also figure \ref{exgts}.

It is clear that each connected component of 
a 
toric graph is a GT graph.
Therefore we define
the partition function
of the toric graph by the product of the partition functions  
of its connected components.

Let us summarize the
information on $X$ read from
the toric graph  $\Gamma$:
\begin{enumerate}
\item
$v\in V_3(\Gamma)$ represents
a torus fixed point $p_v$.
\item
$e\in E_3(\Gamma)$ represents
a curve $C_e\cong \mathbb{P}^1$.
If the two endpoints of $e\in E_3(\Gamma)$ is $v,v'$,
then $p_{v},p_{v'}$ are two torus fixed points in $C_e$.
The framing $n_f$ of $f=(v,e)$ represents the 
degrees of 
the normal bundle:
$N C_e\cong {\mathcal O}_{\mathbb{P}^1}(n_f-1)
\oplus 
{\mathcal O}_{\mathbb{P}^1}(-n_f-1)$.
\end{enumerate}
%
\subsection{Geometric Meaning of Free Energy}%
%
Let $X$ be a TCY threefold.
Roughly speaking,
the (0-pointed) Gromov--Witten invariant
$N_{\beta}^g(X)$
is the number obtained by the integration of $1$ over
the moduli of the 0-pointed stable maps
from a curve of genus $g$ 
whose image belong to the homology class
$\beta\in H_2^{cpt}(X,\mathbb{Z})$.
(See \cite{lllz} for the precise definition.)
We define
the generating function 
with a fixed homology class 
$\beta$:
\begin{equation}\notag
{\mathcal F}_{\beta}(X)=\sum_{g\geq 0}g_s^{2g-2}N_{\beta}^{g}(X).
\end{equation}
In this article, 
we use the symbol $g_s$ as the genus expansion parameter.

Let $\Gamma$ be the toric graph of $X$
and ${\mathcal F}^{\Gamma}_{\Vec{d}}(q)$
be the free energy in a degree $\Vec{d}$.
Note that any 
degree $\Vec{d}$ determines a homology class
with the compact support,
\begin{equation}\notag
[\Vec{d}\cdot \Vec{C}]=\sum_{e\in E_3(\Gamma)}d_e[C_e].
\end{equation}
The proposal of \cite{akmv} (proposition 7.4 \cite{lllz})
is that the 
generating function ${\mathcal F}_{\beta}(X)$ is equal to the
sum of the 
free energy in degrees $\Vec{d}$ such that 
$\Vec{d}\cdot \Vec{C}=\beta$,
under the identification $q=e^{\sqrt{-1}g_s}$:
\begin{equation}\label{betavsd}
{\mathcal F}_{\beta}(X)=\sum_{[\vec{d}.\vec{C}]=\beta}
{\mathcal F}_{\Vec{d}}^{\Gamma}(e^{\sqrt{-1}g_s}).
\end{equation}
Actually, 
each 
${\mathcal F}_{\Vec{d}}^{\Gamma}(q)$
has the meaning in
the localization calculation:
it
is the contribution from the fixed point loci
in the moduli stacks of stable maps
whose image curves are $\Vec{d}\cdot\Vec{C}$.
%
\subsection{Gopakumar--Vafa Conjecture}\label{GVC}
Let us 
define
the numbers
$\{n_{\beta}^g(X)\}_{g\geq 0,\beta\in H_2(X;\mathbb{Z})}$ by rewriting
$\{{\mathcal F}_{\beta}(X)\}_{\beta\in H_2(X;\mathbb{Z})}$
in the form below.
\begin{equation}\label{kakinaosi}
{\mathcal F}_{\beta}(X)=\sum_{g\geq 0}\sum_{k;k|\beta}
\frac{n_{\beta/k}^g(X)}{k}
\Big(2\sin\frac{kg_s}{2}\Big)^{2g-2}.
\end{equation}
$n_{\beta}^g(X)$ is called the {\em Gopakumar--Vafa}
invariant.
The Gopakumar--Vafa conjecture states the followings
\cite{gv}.
\begin{enumerate}
\item
$n_{\beta}^g(X)\in\mathbb{Z}$ and
$n_{\beta}^g(X)=0$ for every fixed $\beta$ and $g\gg1$.
\item
\label{bps}
Moreover,
$n_{\beta}^g(X)$ is equal to the 
number of certain BPS states in M-theory
(see \cite{hst} for a mathematical formulation).
\end{enumerate}

The first part of the conjecture follows from
corollary \ref{mainthm}
since
the GV invariant is written as
\begin{equation}\notag
n_{\beta}^g(X)=\sum_{\Vec{d};[\Vec{d}.\Vec{C}]=\beta}
n_{\Vec{d}}^g(\Gamma),
\end{equation}
by
$(\ref{fn})$,$(\ref{betavsd})$ and $(\ref{kakinaosi})$.
%
%

\section{Proof of Proposition \ref{integrality}}  \label{PI}
In this section, 
we give a proof of proposition 
\ref{integrality}.
%
\subsection{Outline of Proof}
%
The proof proceeds as follows.
Firstly, we
take the logarithm of the partition function
using the Taylor expansion. 
For a degree $\Vec{d}$,
we define 
\begin{equation}\notag
\begin{split}
D(\Vec{d})=\{\Vec{\delta}:\text{ degree}| 
\delta_e\leq d_e \text{ for all } e\in E_3(\Gamma)\}.
\end{split}
\end{equation}
This is just the set of degrees 
smaller than or equal to $\Vec{d}$.
Consider assigning 
a nonnegative integer to each element of $D(\Vec{d})$; 
in other words, consider a set of nonnegative integers
\begin{equation}\notag
n=\{n_{\Vec{\delta}}\in\mathbb{Z}_{\geq 0}
\,|\,
\Vec{\delta}\in D(\Vec{d})\}.
\end{equation}
We call
such a set a {\em multiplicity} in $\Vec{d}$
if it satisfies
\begin{equation}\notag
n\cdot \Vec{d}:=\sum_{\Vec{\delta}\in D(\Vec{d})}
n_{\Vec{\delta}} \, \Vec{\delta}=\Vec{d}.
\end{equation}
With these notations,
the free energy in degree $\Vec{d}$ is written as
\begin{equation}\notag
{\mathcal F}_{\Vec{d}}^{\Gamma}(q)
=\sum_{\begin{subarray}{c}
n;\\
\text{multiplicity}\\
\text{in $\Vec{d}$}  
\end{subarray}}
\frac{|{n}|!}{\prod_{\delta\in D(\Vec{d})} n_{\Vec{\delta}}!}
\frac{(-1)^{|n|-1}}{|n|}
\prod_{\Vec{\delta}\in D(\Vec{d})} 
\big(
{\mathcal Z}_{\Vec{\delta}}^{\Gamma}(q)\big)
^{n_{\Vec{\delta}}}.
\end{equation}
We further rewrite it. 
Let $d_0={\rm gcd}(\Vec{d})$.
\begin{equation}\notag
{\mathcal F}_{\Vec{d}}^{\Gamma}(q)
=
\sum_{k;k|d_0}
\sum_{\begin{subarray}{c}
{n};
\text{multiplicity}\\
\text{in $\Vec{d}/k$},\\
{\rm gcd}(n)=1
\end{subarray}}
\frac{(k|n|)!}{\prod_{\Vec{\delta}\in D(\Vec{d/k})} (kn_{\Vec{\delta}})!}
\frac{(-1)^{k|n|-1}}{k|n|}
\Bigg(\prod_{\Vec{\delta}\in D(\Vec{d}/k)} 
\big(
{\mathcal Z}_{\Vec{\delta}}^{\Gamma}(q)
\big)^{n_{\Vec{\delta}}}
\Bigg)^{k}.
\end{equation}
Then 
\begin{equation}\notag
\begin{split}
G_{\Vec{d}}^{\Gamma}(q)
%
&=
\sum_{k;k|d_0}
\sum_{\begin{subarray}{c}
n;
\text{multiplicity}\\
\text{in $\Vec{d}/k$},\\
{\rm gcd}(n)=1
\end{subarray}}
\frac{1}{k|n|}
\\
&\times
\Bigg[
\sum_{k';k'|k}
\mu\Big(\frac{k}{k'}\Big)
\frac{(k'|n|)!}{\prod_{\Vec{\delta}\in D(\Vec{d}/k)} (k'n_{\Vec{\delta}})!}
\frac{(-1)^{k'|n|-1}}{k'|n|}
\Bigg(\prod_{\Vec{\delta}\in D(\Vec{d}/k)} 
\big(
{\mathcal Z}_{\Vec{\delta}}^{\Gamma}
(q^{k/k'})
\big)^{n_{\Vec{\delta}}}
\Bigg)^{k'}
\Bigg].
\end{split}
\end{equation}
Each summand turns out to be an element of ${\mathcal L}[t]$
by the next lemmas.
\begin{lemma}\label{zinL}%
For any degree $\Vec{d}$,
\begin{equation}\notag
{\mathcal Z}_{\Vec{d}}^{\Gamma}(q)\in{\mathcal L}[t].
\end{equation}
\end{lemma}
%
\begin{lemma}\label{arith}
Let ${n}=(n_1,\ldots, n_l)$ be the set of nonnegative 
integers 
such that ${\rm gcd}({n})=1$.
For 
$R(t)\in {\mathcal L}[t]$,
$k\in \mathbb{N}$ 
and $n$,
\begin{equation}\notag
\begin{split}
\frac{1}{k|n|}
\sum_{k';k'|k}
\mu\Big(\frac{k}{k'}\Big)
\frac{(k'|n|)!}{(k'n_1)!\cdots(k'n_l)!}
\frac{(-1)^{k'|n|}}{k'|n|}
R(t_{k/k'})^{k'}
\in {\mathcal L}[t].
\end{split}
\end{equation}
\end{lemma}
The proofs of lemmas \ref{zinL} and \ref{arith}  
are given in subsection \ref{plz} and
appendix \ref{proofarith}, respectively.

Thus $G_{\Vec{d}}^{\Gamma}(q)\in{\mathcal L}[t]$
and
proposition \ref{integrality} is proved.

\subsection{Proof of Lemma \ref{zinL}}\label{plz}
 In this subsection,
we give
a proof of lemma \ref{zinL}.
The main point is in showing  that
${\mathcal Z}_{\Vec{d}}^{\Gamma}(q)$,
which is a priori a function in $q^{\frac{1}{2}}$,
is actually a function in $t$.
%
%
%
We use two key facts here.
Let $\mathbb{Z}_0[t]$ be the ring of monic polynomials
and 
let $\mathbb{Z}^+[q,q^{-1}]$
be the subring 
of the ring of Laurent polynomials in $q$
whose elements are
symmetric with respect to $q,q^{-1}$.
The one fact is that \cite{bp}
\begin{equation}\notag
t_k:=[k]^2\in\mathbb{Z}_0[t]
\qquad
(k\in\mathbb{N}).
\end{equation}
The other is that (see \cite{ko}, lemma 6.2)
\begin{equation}\notag
{\mathbb{Z}}[t]\cong \mathbb{Z}^+[q,q^{-1}].
\end{equation}
%

We first state preliminary lemmas.
\renewcommand{\theenumi}{\roman{enumi}}
\renewcommand{\labelenumi}{(\theenumi)}
\begin{lemma}\label{prem}
\begin{enumerate}
\item \label{ichi}
$h_i(q^{\rho})$ is written in the form
\begin{equation}\notag
h_i(q^{\rho})=q^{i/2}\frac{f_2(q)}{f_1(t)}
\end{equation}
with $f_2(q)\in \mathbb{Z}[q,q^{-1}]$
and $f_1(q)\in \mathbb{Z}_0[t]$.
\item \label{ni}
$e_i(q^{\rho})=(-1)^i h_i(q^{\rho})|_{q\to q^{-1}}$.
\item \label{san}
$h_i(q^{\lambda+\rho})$ is written in the form
\begin{equation}\notag
h_i(q^{\lambda+\rho})=q^{i/2}\frac{f_2^{\lambda}(q)}{f_1^{\lambda}(t)}
\end{equation}
with  $f_2^{\lambda}(q)\in \mathbb{Z}[q,q^{-1}]$
and $f_1^{\lambda}(t)\in \mathbb{Z}_0[t]$.
\item 
\label{yon}
$e_i(q^{\lambda+\rho})
=(-1)^i h_i(q^{\lambda^t+\rho})|_{q\to q^{-1}}$.
\item\label{go}
$s_{\mu/\nu}(q^{\lambda+\rho})$ is written in the following
form:
\begin{equation}\notag
s_{\mu/\nu}(q^{\lambda+\rho})
=q^{\frac{|\mu|-|\nu|}{2}}\frac{s_2^{\lambda,\mu,\nu}(q)}
{s_1^{\lambda,\mu,\nu}(t)}
\end{equation}
with 
$s_2^{\lambda,\mu,\nu}(q)\in\mathbb{Z}^+[q,q^{-1}]$
and $s_1^{\lambda,\mu,\nu}(t)\in\mathbb{Z}_0[t]$.
\item 
\label{roku}
$s_{\mu^t/\nu^t}(q^{\lambda+\rho})=
(-1)^{|\mu|-|\nu|}s_{\mu/\nu}(q^{\lambda+\rho})|_{q\to q^{-1}}
$.
\item \label{nana}
The three point function is written in the following form:
\begin{equation}\notag
C_{\lambda^1,\lambda^2,\lambda^3}(q)
=q^{\frac{|\lambda^1|+|\lambda|^2+|\lambda^3|}{2}}
\frac{c_2^{\lambda^1,\lambda^2,\lambda^3}(q)}
{c_1^{\lambda^1,\lambda^2,\lambda^3}(t)}
\end{equation}
where
$c_2^{\lambda^1,\lambda^2,\lambda^3}(q)\in \mathbb{Z}[q,q^{-1}]$ 
and
$c_1^{\lambda^1,\lambda^2,\lambda^3}(t)\in\mathbb{Z}_0[t]$.
\item \label{hachi}
\begin{equation}\notag
C_{\lambda^{1t},\lambda^{2t},\lambda^{3t}}(q)=
(-1)^{|\lambda^1|+|\lambda^2|+|\lambda^3|}
C_{\lambda^{1},\lambda^{2},\lambda^{3}}(q^{-1}).
\end{equation}
\end{enumerate}
\end{lemma}
\renewcommand{\theenumi}{\arabic{enumi}}
\renewcommand{\labelenumi}{\theenumi.}

\begin{proof}%
(\ref{ichi}).
Recall the expression (\ref{eh1}).
If we multiply both the denominator and the numerator 
by
$[1]\ldots[i]$, we obtain
\begin{equation}\notag
h_i(q^{\rho})=q^{\frac{i}{2}}
\frac{q^{\frac{i(i-3)}{4}}[1]\cdots[i]}{t_1\cdots t_i}.
\end{equation}
This proves (\ref{ichi}).

(\ref{ni}) follows from (\ref{eh1}).

(\ref{san}) follows from (\ref{ichi}) and
the generating function (\ref{eh2}).
 
(\ref{yon}) follows from
(\ref{eh2}) and
the identity:
\begin{equation}\notag
\prod_{i=1}^{l(\lambda)}
\frac{1+q^{\lambda_i-i+\frac{1}{2}}z}
{1+q^{-i+\frac{1}{2}}z} 
=
\prod_{j=1}^{l(\lambda^t)}
\frac{1+q^{j-\frac{1}{2}}z}
{1+q^{-\lambda_j^t+j-\frac{1}{2}}z}.
\end{equation}
(This identity can be proved by showing that
the LHS is equal to 
$\prod_{i=1}^{r(\lambda)}
(1+{q^{\lambda_i-i+\frac{1}{2}}z})/
(1+q^{-(\lambda_i^t-i+\frac{1}{2})z})
$
where
$r(\lambda)$ denotes
the number of diagonal boxes in the Young diagram of $\lambda$.)

(\ref{go})
follows 
from (\ref{san}) and (\ref{she}):
\begin{equation}\notag
s_{\mu/\nu}(q^{\lambda+\rho})
=
\det\big(h_{\mu_i-\nu_j-i+j}(q^{\lambda+\rho})\big)_{i,j}
=
q^{\frac{|\mu|-|\nu|}{2}}
\det\big(q^{\frac{-\mu_i+\nu_j+i-j}{2}}
h_{\mu_i-\nu_j-i+j}(q^{\lambda+\rho})\big)_{i,j}.
\end{equation}

(\ref{roku})
follows from (\ref{yon}) and (\ref{she})
\begin{equation}\notag
\begin{split}
s_{\mu^t/\nu^t}(q^{\lambda+\rho})
&=
\det\big(e_{\mu_i-\nu_j-i+j}(q^{\lambda+\rho})\big)_{i,j}
\\
&=
(-1)^{|\mu|-|\nu|}
\det\big(
h_{\mu_i-\nu_j-i+j}(q^{\lambda^t+\rho})\big)_{i,j}|_{q\to q^{-1}}
\qquad (\because \text{(\ref{yon})})
\\
&=(-1)^{|\mu|-|\nu|}s_{\mu/\nu}(q^{\lambda^t+\rho})|_{q\to q^{-1}}.
\end{split}
\end{equation}

(\ref{nana}) and (\ref{hachi}) 
follow from (\ref{go}) and 
(\ref{roku}), respectively.
\end{proof}

Now we prove
lemma \ref{zinL}.
Let $\Vec{\lambda}=(\lambda_e)_{e\in E_3(\Gamma)}$ 
be a $\Gamma$-partition.
By (\ref{nana}),
$Y_{\Vec{\lambda}}(q)$ (defined in (\ref{defy}))
is written in the form
\begin{equation}\notag
Y_{\Vec{\lambda}}(q)=
\frac{Y_2^{\Vec{\lambda}}(q)}
{Y_1^{\Vec{\lambda}}(t)}
\end{equation}
with 
$Y_2^{\Vec{\lambda}}(q)\in \mathbb{Z}[q,q^{-1}]$
and 
$Y_1^{\Vec{\lambda}}(t)\in \mathbb{Z}_0[t]$.
Moreover, by (\ref{hachi}), it holds 
that
\begin{equation}\notag
Y_{\Vec{\lambda}^t}(q)=
\frac{Y_2^{\Vec{\lambda}}(q^{-1})}
{Y_1^{\Vec{\lambda}}(t)}
=Y_{\Vec{\lambda}}(q^{-1}),
\end{equation}
where $\Vec{\lambda}^t =({\lambda_e}^t)_{e\in E_3(\Gamma)}$.
Therefore 
\begin{equation}\notag
Y_{\Vec{\lambda}}(q)+Y_{{\Vec{\lambda}}^t}(q)
\in {\mathcal L}[t] 
\quad (\Vec{\lambda}\neq \Vec{\lambda}^t)
,\qquad
Y_{\Vec{\lambda}}(q)\in {\mathcal L}[t]
\quad (\Vec{\lambda}= \Vec{\lambda}^t).
\end{equation}
Thus
\begin{equation}\notag
{\mathcal Z}_{\Vec{d}}^{\Gamma}(q)=
\frac{1}{2}
\sum_{\begin{subarray}{c}
\Vec{\lambda}; \text{ $\Gamma$-partition}\\
\text{ of degree } \Vec{d},\\
\Vec{\lambda}\neq \Vec{\lambda}^t
\end{subarray}
}
\Big(Y_{\Vec{\lambda}}(q)+Y_{\Vec{\lambda}^t}(q)
\Big)+
\sum_{\begin{subarray}{c}
\Vec{\lambda}; \text{ $\Gamma$-partition}\\
\text{ of degree }\Vec{d},\\
\Vec{\lambda}=\Vec{\lambda^t}
\end{subarray}
}
Y_{\Vec{\lambda}}(q)\in {\mathcal L}[t].
\end{equation}
Note that
the prefactor $1/2$ 
does not matter because
the same term appears twice
if $\Vec{\lambda}\neq \Vec{\lambda}^t$.
The proof of  lemma \ref{zinL} is finished.

\section{Proof of Proposition \ref{polynomiality}}\label{PV}
In this section, we give a proof of 
proposition \ref{polynomiality}.
We first 
rewrite the three point function 
and the partition function
in the operator formalism (subsections \ref{tpf} and \ref{pf}).
Then we express the partition function as the sum of 
certain quantities - {\em combined amplitude} - of
not necessarily connected graphs
(subsection \ref{ge}).
By using the exponential formula,
we
obtain the free energy as the sum over connected graphs
(subsection \ref{feen}).
Then we show that 
the proposition follows from the property of 
the combined amplitudes of the connected graphs
(subsection \ref{ppp}).
This proof is  almost the same as \cite{ko},
where the same proposition was proved for
the middle graph in figure \ref{exgts}.

There are, however, two technical difficulties in generalization.
They occur
when writing the partition function in the operator formalism
due to the existence of
trivalent vertices whose three incident edges are in $E_3(\Gamma)$.
The one difficulty is
how to incorporate the variables 
such as $q^{\lambda+\rho}$.
%
It is overcome by the fact that
the $i$-th power sum of $q^{\lambda+\rho}$ is equal to
the matrix element of ${\mathcal E}_0(i)$ with respect to the
state $|v_{\lambda}\rangle$ (subsection \ref{ps}).
The other is how to perform the summation 
when the states
$|v_{\lambda}\rangle$ and $|v_{\lambda^t}\rangle$
appear simultaneously.
It is  solved by introducing the operator $R$
that transforms $|v_{\lambda}\rangle$ to 
$|v_{\lambda^t}\rangle$ (subsection \ref{trop}).

The expression of the 
three point function $C_{\lambda^1,\lambda^2,\lambda^3}$
thus obtained
does not posses the cyclic symmetry
with respect to three partitions 
$\lambda^1,\lambda^2,\lambda^3$.
Therefore it becomes necessary to specify the order
of three flags around every trivalent vertex.
For this purpose we will introduce the notion of 
the {\em flag-order}
(subsection \ref{flagorder}).

We omit the explanation of 
the operator formalism.
Please see \cite{ko}, section {2.1}.

%

\subsection{Technical Preliminary}\label{tec}
This subsection is devoted to the solution to 
the two technical problems mentioned  previously.

%

\subsubsection{Power Sum}\label{ps}
We express the power sum
functions of the variables $q^{\lambda+\rho}$
as a matrix element.

For a sequence of variables $x=(x_1,x_2,\ldots)$,
the $i$-th power sum function is defined by 
$p_i(x)=\sum_{j\geq 1}{x_j}^i$.
The power sum function associated to a partition $\nu$
is defined by $p_{\nu}(x)=\prod_{i=1}^{l(\nu)}p_{\nu_i}(x)$.

Consider the variable 
$q^{\lambda+\rho}=(q^{\lambda_i-i+\frac{1}{2}})_{i\geq 1}$
associated to a partition $\lambda$.
The $i$-th power sum function is equal to
\begin{equation}\label{pi}
\begin{split}
p_{i}(q^{\lambda+\rho})&=
\underbrace{
\sum_{j=1}^{l(\lambda)}
\big(
q^{i(\lambda_j-j+\frac{1}{2})}-q^{i(-j+\frac{1}{2})}
\big)
}_{\star}
+\frac{1}{[i]} \qquad (i\in \mathbb{Z}\setminus\{0\}).
\end{split}
\end{equation}
It turns out to be written 
as the matrix element of 
the operator ${\mathcal E}_0(i)$,
\begin{equation}\notag
{\mathcal E}_0(i)=\sum_{k\in\halfintegers}q^{ik}E_{k,k}+\frac{1}{[i]}
 \qquad (i\in \mathbb{Z}\setminus\{0\}).
\end{equation}
\begin{lemma}\label{powersum}
\begin{equation}\notag
p_i(q^{\lambda+\rho})=\langle v_{\lambda}|
{\mathcal E}_0(i)|v_{\lambda}\rangle
=
-\langle v_{\lambda^t}|{\mathcal E}_0(-i)|v_{\lambda^t}\rangle.
\end{equation}
\end{lemma}

%
The lemma implies that
\begin{equation}\notag
p_{\nu}(q^{\lambda+\rho})=
\langle v_{\lambda}|{\mathcal E}_0(\nu)|v_{\lambda}\rangle,
\qquad
p_{\nu}(q^{\lambda^t+\rho})=(-1)^{l(\nu)}
\langle v_{\lambda}|{\mathcal E}_0(-\nu)|v_{\lambda}\rangle,
\end{equation}
where
\begin{equation}\label{Ezero}
{\mathcal E}_{0}(\pm\nu):=
{\mathcal E}_{0}(\pm \nu_1)\cdots {\mathcal E}_0(\pm \nu_{l(\nu)}).
\end{equation}
\begin{proof}
%
%
Note that
the state $|v_{\lambda}\rangle$ is 
written as
\begin{equation}\notag
|v_{\lambda}\rangle
=
\psi_{\lambda_1-\frac{1}{2}}\ldots
\psi_{\lambda_{l(\lambda)}-l(\lambda)+\frac{1}{2}}
\psi^*_{-l(\lambda)+\frac{1}{2}}
\ldots\psi^*_{-\frac{1}{2}}
|0\rangle.
\end{equation}
Therefore the action of ${\mathcal E}_0(i)$ on $|v_{\lambda}\rangle$
follows from 
the commutation relations:
\begin{equation}\notag
[{\mathcal E}_0(i),\psi_k]=q^{ik}\psi_k,\qquad
[{\mathcal E}_0(i),\psi_k^*]=-q^{ik}\psi_k^*.
\end{equation}
Thus we have
\begin{equation}\notag
{\mathcal E}_0(i)|v_{\lambda}\rangle
=p_i(q^{\lambda+\rho})|v_{\lambda}\rangle.
\end{equation}
So if we take the pairing with
$\langle v_{\lambda}|$, we obtain
the first line.

To show the second,
we need the next identities.
The term $(\star)$ in (\ref{pi}) is written 
in three ways:
\begin{equation}\notag
\begin{split}
(\star)&=
\sum_{j=1}^{l(\lambda)}
\big(
q^{i(\lambda_j-j+\frac{1}{2})}-q^{i(-j+\frac{1}{2})}
\big)
\\
&=
\sum_{j=1}^{r(\lambda)}
\big(
q^{i(\lambda_j-j+\frac{1}{2})}-q^{-i(\lambda_j^t-j+\frac{1}{2})}
\big)
\\
&=-
\sum_{j=1}^{l(\lambda^t)}
\big(
q^{-i(\lambda_j^t-j+\frac{1}{2})}-q^{-i(-j+\frac{1}{2})}
\big).
\end{split}
\end{equation}
In the second line, $r(\lambda)$ denotes
the number of diagonal boxes in the Young diagram of $\lambda$.
This implies
\begin{equation}\notag
p_i(q^{\lambda+\rho})=-p_i(q^{\lambda^t+\rho})|_{q\to q^{-1}}
=-\langle v_{\lambda^t}|{\mathcal E}_0(-i)|v_{\lambda^t}\rangle.
\end{equation}
\end{proof}

\subsubsection{Transposition Operator}\label{trop}
We define the
{\em transposition operator}
$R$ as follows.
The actions of $R$ on the 
charge zero subspace
$\Lambda^{\frac{\infty}{2}}_0V$
and its dual space are defined by
\begin{equation}\notag 
R|v_{\lambda}\rangle =R|v_{\lambda}^t\rangle,\qquad
\langle v_{\lambda}|R=\langle v_{\lambda}^t|.
\end{equation}
The actions of $R$ on the fermions
are determined by compatibility%
\footnote{
One could check the compatibility by using the following
expressions:
\begin{equation}\notag
|v_{\lambda}\rangle =\prod_{i=1}^{r(\lambda)}
\big((-1)^{b_i-\frac{1}{2}}\psi_{a_i}\psi^*_{-b_i}\big)|0\rangle,
\qquad
\langle v_{\lambda}|=
\langle 0|\prod_{i=1}^{r(\lambda)}
\big((-1)^{b_i-\frac{1}{2}}\psi_{b_i}\psi^*_{a_i}\big)
\end{equation}
where $r(\lambda)=\#\text{ (diagonal boxes in $\lambda$)}$,
$a_i=\lambda_i-i+\frac{1}{2}$, $b_i=\lambda^{t}_i-i+\frac{1}{2}$.}%
:
\begin{equation}\notag
R\psi_k R^{-1}=(-1)^{k-\frac{1}{2}}\psi_{-k}^*,\qquad
R\psi_k^* R^{-1}=(-1)^{-k+\frac{1}{2}}\psi_{-k}
\qquad(k\in\halfintegers).
\end{equation}
Therefore the  actions on other operators are as follows:
\begin{equation}\label{ar}
\begin{split}
&R E_{i,j}R^{-1}=(-1)^{i-j+1}E_{-j,-i},
\\
&R{\mathcal E}_c(n)R^{-1}=
(-1)^{c+1}{\mathcal E}_c(-n)\qquad (c,n)\neq (0,0),
\\
&R \alpha_m R^{-1}=(-1)^{m+1}\alpha_{m} \qquad
(m\in\mathbb{Z},m\neq 0),
\\
& R {\mathcal F}_2 R^{-1}=
-{\mathcal F}_2, 
\qquad
RHR^{-1}=H.
\end{split}
\end{equation}
The action on a bosonic state $|\mu\rangle$ is obtained from 
the third of (\ref{ar}):
\begin{equation}\label{aar}
R|\mu\rangle=(-1)^{l(\mu)+|\mu|}|\mu\rangle.
\end{equation}

\subsection{Three Point Function}\label{tpf}
%
The three point function $C_{\lambda^1,\lambda^2,\lambda^3}(q)$ 
is written in the operator
formalism as follows.
\begin{lemma}\label{threepoint}
\begin{equation}\notag
\begin{split}
C_{\lambda^1,\lambda^2,\lambda^3}(q)
&=
\sum_{\begin{subarray}{c}\mu\in\mathcal P\\
|\mu|\leq |\lambda^1|,|\lambda^3|
\end{subarray}}
\sum_{\begin{subarray}{c}
\nu^1,\nu^2,\nu^3\in{\mathcal P};\\
|\nu^2|=|\lambda^2|,\\
|\nu^1|=|\lambda^1|-|\mu|,\\
|\nu^3|=|\lambda^3|-|\mu|
\end{subarray}}
\frac{(-1)^{l(\nu^1)}}
{z_{\mu}z_{\nu^1}z_{\nu^2}z_{\nu^3}[\nu^2]}
\\&\times
\langle v_{\lambda^1}|\nu^1\cup\mu\rangle
\langle v_{\lambda^{3t}}|q^{-{\mathcal F}_2}|\nu^3\cup\mu\rangle
\langle v_{\lambda^2}|{\mathcal E}_0(-\nu^1){\mathcal E}_0(\nu^3)
|\nu^2\rangle.
\end{split}
\end{equation}
\end{lemma}
Note that the cyclic symmetry 
with respect to $\lambda^1,\lambda^2,\lambda^3$
is not manifest in  this expression.
%
\begin{proof}%
The skew-Schur function 
in the variables $x=(x_1,x_2,\ldots)$ is written as
\begin{equation}\notag
s_{\mu/\eta}(x)=
\sum_{\begin{subarray}{c}\mu';|\mu'|=|\mu|-|\eta|,\\
                         \eta';|\eta'|=|\eta|
      \end{subarray}}
\frac{p_{\mu'}(x)}{z_{\mu'}z_{\eta'}}
\langle v_{\mu}|\mu'\cup\eta'\rangle
\langle \eta'|v_{\eta}\rangle.
\end{equation}
Therefore,
\begin{equation}\notag
s_{\lambda^2}(q^{\rho})=
\sum_{\begin{subarray}{c}\nu^2\in{\mathcal P};
|\nu^2|=|\lambda^2|\end{subarray}}
\frac{1}{z_{\nu^2}}
\langle v_{\lambda^2}|\nu^2\rangle.
\end{equation}
And
\begin{equation}\notag
\begin{split}
&\sum_{\eta} s_{\lambda^1/\eta}(q^{\lambda^{2t}+\rho})
s_{\lambda^{3t}/\eta}(q^{\lambda^2+\rho})\\
&=
\sum_{d=0}^{\min\{|\lambda^1|,|\lambda^3|\}}
\sum_{\eta\in{\mathcal P}_d}
s_{\lambda^1/\eta}(q^{\lambda^{2t}+\rho})
s_{\lambda^{3t}/\eta}(q^{\lambda^2+\rho})
\\
&=\sum_{d=0}^{\min\{|\lambda^1|,|\lambda^3|\}}
\sum_{\begin{subarray}{c}
\nu^1;|\nu^1|=|\lambda^1|-d,\\
\nu^3;|\nu^3|=|\lambda^3|-d,\\
\mu,\mu'\in{\mathcal P}_d
\end{subarray}}
\frac{
p_{\nu^1}(q^{\lambda^{2t}+\rho})
p_{\nu^3}(q^{\lambda^2+\rho})
}
{z_{\nu^1}z_{\nu^3}z_{\mu}z_{\mu'}}
\langle v_{\lambda^1}|\nu^1\cup\mu\rangle
\langle v_{\lambda^{3t}}|\nu^3\cup\mu'\rangle
\underbrace{
\sum_{\eta\in{\mathcal P}_d}
\langle\mu|v_{\eta} \rangle
\langle v_{\eta}|\mu' \rangle
}_{=\langle\mu|\mu'\rangle=z_{\mu}\delta_{\mu,\mu'}}.
\end{split}
\end{equation}
By lemma \ref{powersum},
this is equal to
\begin{equation}\notag
\begin{split}
\sum_{\nu^1,\nu^3,\mu}
\frac{(-1)^{l(\nu^1)}}
{z_{\nu^1}z_{\nu^2}z_{\mu}}
\langle v_{\lambda^1}|\nu^1\cup\mu\rangle
\langle v_{\lambda^{3t}}|\nu^3\cup\mu\rangle
\langle v_{\lambda^2}|{\mathcal E}_0(-\nu^1){\mathcal E}_0(\nu^3)
|v_{\lambda^2}\rangle.
\end{split}
\end{equation}
The factor $q^{\frac{\kappa(\lambda_3)}{2}}$ is written
as follows:
\begin{equation}\notag
q^{\frac{\kappa(\lambda^3)}{2}}=
\langle v_{\lambda^{3t}}|q^{-{\mathcal F}_2}|v_{\lambda^{3t}}\rangle.
\end{equation}
Combining the above expressions, 
we obtain the lemma.
\end{proof}
\subsection{Flag-order}\label{flagorder}
Since the expression in 
lemma \ref{threepoint}
is not cyclic symmetric, 
we have to
specify 
a counterclockwise order of three flags
for
every trivalent vertex $v$.
For this reason, we  introduce the notion of 
{\em flag order} of a GT graph $\Gamma$.

Let 
$F(v)$ be the set of three flags incident on a
trivalent vertex $v\in V_3(\Gamma)$
and
$F_3'(\Gamma)$ be the set of flags whose vertices are 
trivalent:
\begin{equation}\notag
F_3'(\Gamma)=\bigcup_{v\in V_3(\Gamma)}F(v).
\end{equation}
(Note that
$F_3(\Gamma)\subset F_3'(\Gamma)$.)
A {\em flag-order} of a GT graph $\Gamma$ is a map 
$\iota$ from 
$F_3'(\Gamma)$ to $\{1,2,3\}$
satisfying the following conditions: 
for every $v\in V_3(\Gamma)$,
\begin{enumerate}
\item
$\iota_v:F(v)\to \{1,2,3\}$ is one-to-one;
\item
the disposition of three flags
$\iota_v^{-1}(1),\iota_v^{-1}(2),\iota_v^{-1}(3)$
is counterclockwise,
\end{enumerate}
where $\iota_v$ is the restriction of $\iota$ to $F(v)$.
\newcommand{\myf}{\mathfrak{f}}
For convenience of writing, we set
$\myf_i(v):=\iota_v^{-1}(i)$ $(i=1,2,3)$
for $v\in V_3(\Gamma)$.
We also set
$\myf_i(f):=\myf_i(v)$ for $f\in F_3(\Gamma)$
where $v$ is the vertex of $f$.

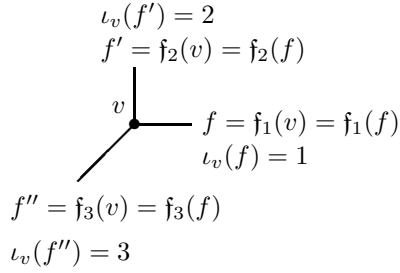
\newcommand{\flagorder}{
\unitlength .15cm
\begin{picture}(-10,-15)(8,10)
\thicklines
\put(0,0){\line(1,0){5}}
\put(0,0){\line(0,1){5}}
\put(0,0){\line(-1,-1){5}}
\put(0,0){\circle*{1}}
\put(-2,1){$v$}
\put(6,-.5){$f=\myf_1(v)=\myf_1(f)$}
\put(-3,6){$f'=\myf_2(v)=\myf_2(f)$}
\put(-11,-8){$f''=\myf_3(v)=\myf_3(f)$}
\put(6,-3.5){$\iota_v(f)=1$}
\put(-3,9){$\iota_v(f')=2$}
\put(-11,-12){$\iota_v(f'')=3$}
\end{picture}
}
\begin{equation}\notag
\raisebox{2.5cm}{
\flagorder}
\end{equation}

\vspace*{1cm}

A flag-order of a GT graph $\Gamma$ is not unique;
there are 
$3^{\#V_3(\Gamma)}$ flag-orders.

In the rest of section \ref{PV}, 
we fix one GT graph $\Gamma$ and 
one flag-order $\iota$.

\subsection{Partition Function}\label{pf}
%
The goal of this subsection is to
rewrite the partition function 
${\mathcal Z}_{\vec{d}}^{\Gamma}(q)$ in the operator formalism.
To state the result, we 
introduce the following notations.

\begin{itemize}
\item
We use the symbol $\Vec{\mu}$ and $\vec{\nu}$
for tuples of partitions
$\vec{\mu}=(\mu^v)_{v\in V_3(\Gamma)}$ 
and $\vec{\nu}=(\nu^f)_{f\in F_3(\Gamma)}$.
A {\em $\Gamma$-set of degree $\Vec{d}$}
is a pair $(\vec{\mu},\Vec{\nu})$ satisfying
the following conditions:
\begin{enumerate}
\item
for $f=(v,e)\in F_3(\Gamma)$,
$|\mu^v|+|\nu^f|=d_e$ if  $\iota(f)=1,3$;
\item
for $f=(v,e)\in F_3(\Gamma)$,
$|\nu^f|=d_e$ if $\iota(f)=2$;
\item
for $v\in V_3(\Gamma)$,  
$|\mu^v|=0$
if $\myf_1(v)\notin F_3(\Gamma)$ 
or $\myf_3(v)\notin F_3(\Gamma)$.
\end{enumerate}
%
%
\item
Integers $N_e$,
$L_e(\vec{\mu},\vec{\nu})$
and matrix elements $K_e(\vec{\mu},\vec{\nu})$ 
($e\in E_3(\Gamma)$)
are defined 
as in table \ref{tab:NL}.
\end{itemize}
With these notations,
${\mathcal Z}_{\Vec{d}}^{\Gamma}(q)$
is written as follows:
\begin{lemma}\label{zop}
\begin{equation}\label{zopf}
\begin{split}
 {\mathcal Z}_{\Vec{d}}^{\Gamma}(q)
&=
   \sum_{\begin{subarray}{c}
          (\Vec{\mu},\Vec{\nu});\text{ $\Gamma$-set}\\
          \text{of degree $\vec{d}$}
         \end{subarray}}
   \frac{(-1)^{L_1(\vec{\mu},\vec{\nu})+L_2(\vec{\mu},\vec{\nu})}}
{z_{\vec{\mu}}z_{\vec{\nu}}}
\prod_{\begin{subarray}{c}
                  f\in F_3(\Gamma),\\
                  \iota(f)=2
       \end{subarray}}
\frac{1}{[\nu^f]}
\prod_{\begin{subarray}{c}
            f\in F_3(\Gamma);\\
            \iota(f)=1,3,\\
            \myf_2(f)\notin F_3(\Gamma)
       \end{subarray}}
\frac{1}{[\nu^f]}
   \prod_{e\in E_3(\Gamma)} K_e(\Vec{\mu},\Vec{\nu}).
\end{split}
\end{equation}
Here 
$L_1(\vec{\mu},\vec{\nu})$
and $L_2(\vec{\mu},\vec{\nu})$ 
are
\begin{equation}\notag
L_1(\vec{\mu},\vec{\nu})
=
\sum_{e\in E_3(\Gamma)}L_e(\vec{\mu},\vec{\nu})+
\sum_{\begin{subarray}{c}
            f\in F_3(\Gamma);\\
            \iota(f)=1,\\
            \myf_2(f)\notin F_3(\Gamma)
       \end{subarray}}l(\nu^f)
,\qquad
L_2(\vec{\mu},\vec{\nu})=\sum_{e\in E_3(\Gamma)}N_e d_e,
\end{equation}
and 
$z_{\vec{\mu}}=\prod_{v\in V_3(\Gamma)}z_{\mu^v}$,
$z_{\vec{\nu}}=\prod_{f\in F_3(\Gamma)}z_{\nu^f}$.
\end{lemma}

\begin{table}[t]
\begin{equation}\notag
\begin{array}{|c|c|c|c|}\hline
   (\iota(f),\iota(f'))
  &N_e&L_e(\vec{\mu},\vec{\nu})&K_e(\vec{\mu},\vec{\nu})
\\\hline
(1,1)&-n_e & l(\mu)+l(\nu')
     &\langle\mu\cup \nu |q^{N_e{\mathcal F}_2}|\mu'\cup\nu'\rangle
\\
(1,2)&-n_e & l(\mu)
     &\langle \mu\cup \nu|{\mathcal E}_0(-\nu^{\prime 1})
            {\mathcal E}_0(\nu^{\prime 3})
               q^{N_e{\mathcal F}_2}
             |\nu'\rangle
\\
(1,3)&n_e-1& l(\nu)
     & \langle \mu\cup \nu|
               q^{N_e {\mathcal F}_2}
                      |\mu'\cup \nu'\rangle
\\
(2,1)& -n_e & l(\nu)+l(\nu')+l(\nu^1)+l(\nu^3)
     & \langle \nu|
         {\mathcal E}_0(-\nu^{3}){\mathcal E}_0(\nu^{1})
         q^{N_e{\mathcal F}_2} 
              |\mu'\cup \nu'\rangle
\\
(2,2)&-n_e & l(\nu)+l(\nu^1)+l(\nu^3)
     &\langle \nu|
        {\mathcal E}_0(-\nu^3\cup\nu^{\prime 1})
        {\mathcal E}_0(\nu^1 \cup \nu^{\prime 3})
         q^{N_e{\mathcal F}_2} 
              |\nu'\rangle
\\
(2,3)&n_e-1 &0
     &\langle \nu|
      {\mathcal E}_0(-\nu^{1}){\mathcal E}_0(\nu^{3})
         q^{N_e{\mathcal F}_2} 
              |\mu'\cup \nu'\rangle
\\
(3,1)&-n_e-1 & l(\nu')
     & \langle \mu\cup \nu|
         q^{N_e{\mathcal F}_2} 
              |\mu'\cup \nu'\rangle
\\
(3,2)&-n_e-1 & 0
     &\langle \mu\cup\nu|{\mathcal E}_0(-\nu^{\prime 1})
                         {\mathcal E}_0(\nu^{\prime 3})
         q^{N_e{\mathcal F}_2} 
              |\nu'\rangle
\\
(3,3)&n_e &l(\mu)+l(\nu)
     & \langle \mu\cup \nu|
         q^{N_e{\mathcal F}_2} 
              |\mu'\cup \nu'\rangle
\\\hline
\end{array}
\end{equation}
\newcommand{\Ke}{
\unitlength .2cm
\begin{picture}(12,6)(-6,3)
\thicklines
\put(0,0){\line(1,0){5}}
\put(0,0){\line(-3,5){2}}
\put(0,0){\line(-3,-5){2}}
\put(5,0){\line(3,5){2}}
\put(5,0){\line(3,-5){2}}
\put(0,0){\circle*{.5}}
\put(5,0){\circle*{.5}}
\put(0,0){\vector(1,0){3}}
\put(-2,0){$v$}
\put(6,0){$v'$}
\put(2,1){$e$}
\put(-11,0){$\mu^v=\mu$}
\put(10,0){$\mu^{v'}=\mu'$}
\put(-11,-3){$\nu^{(v,e)}=\nu$}
\put(10,-3){$\nu^{(v',e)}=\nu'$}
\end{picture}
}
\begin{equation}\notag
\raisebox{1.3cm}{
\Ke}
\end{equation}

\caption{$N_e$,$L_e(\vec{\mu},\vec{\nu})$ 
and $K_e(\vec{\mu},\vec{\nu})$.
$\mu,\mu',\nu,\nu'$ are as above.
When $\iota((v,e))=2$,
$\nu^1=
\nu^{\myf_1(v)}$
and
when $\iota((v',e))=2$,
$\nu^{\prime 1}=\nu^{\myf_1(v')}$,
$\nu^{\prime 3}=\nu^{\myf_3(v')}$.
} 
\label{tab:NL}
\end{table}

The details of the RHS of (\ref{zopf}) 
depends on the choice of the flag order $\iota$
(cf. example \ref{exof}-(i)).
However,
the rest of the proof of proposition \ref{polynomiality}
proceeds completely in the same way.

Since (\ref{zopf}) is quite involved,
we compute it for the examples of subsection \ref{expart}.
Then we 
describe a proof in the case of example \ref{expart}-3.
The proof for 
a general GT graph is the same and left to the reader.

\begin{example}\label{exof}
(i) Let $\Gamma$ be the GT graph of example \ref{expart}-1.

We take the flag order
$
\iota((v,e))=1,\iota((v',e))=3.
$ 
The value of $\iota$ for other flags are determined by
the counterclockwise condition.
The partition function with respect to this 
flag order is:
\begin{equation}\notag
{\mathcal Z}_d^{\Gamma}(q)=\sum_{\nu,\nu'\in{\mathcal P}_d}
\frac{(-1)^{d(n+1)}}{z_{\nu}z_{\nu'}
[\nu][\nu']}
\langle \nu|q^{(n-1){\mathcal F}_2}|\nu'\rangle.
\end{equation}
If we take another flag order
$
\iota((v,e))=\iota((v',e))=1,
$ 
then we obtain another expression:
\begin{equation}\notag
{\mathcal Z}_d^{\Gamma}(q)=\sum_{\nu,\nu'\in{\mathcal P}_d}
\frac{(-1)^{dn+l(\nu)}}{z_{\nu}z_{\nu'}
[\nu][\nu']}
\langle \nu|q^{-n{\mathcal F}_2}|\nu'\rangle.
\end{equation}
In this exmaple,
one could see that 
the operator formalism expression 
depends on the choice of the flag order.

(ii) Let $\Gamma$ be the GT graph of
example \ref{expart}-2.
We take the flag-order
\begin{equation}\notag
\iota(v_{i-1},e_i)=3,\qquad \iota(v_{i},e_i)=1
\qquad(1\leq i\leq r).
\end{equation}
We write
$\vec{\mu}=(\mu^1,\ldots,\mu^r)$
where 
$\mu^i$ is a partition associated to the vertex $v^i$.
We also write
$\vec{\nu}=(\vec{\alpha},\vec{\beta})$,
$\vec{\alpha}=(\alpha_1,\ldots,\alpha^r)$,
$\vec{\beta}=(\beta^1,\ldots,\beta^r)$
where
$\alpha^i$ and $\beta^i$ are partitions associated to
flags $(v_{i-1},e_i)$ and  $(v_{i},e_i)$ 
(here $v^{-1}=v^r$ is assumed).
The condition that a 
$\Gamma$-set $(\vec{\mu},\vec{\nu})$ is
of degree $\vec{d}=(d_1,\ldots,d_r)$
is
\begin{equation}\notag
|\mu^{i-1}|+|\alpha^i|=
|\mu^{i}|+|\beta^i|=d_i \qquad(1\leq i\leq r),
\end{equation}
where $\mu^{-1}=\mu^r$.
Then
\begin{equation}\notag
{\mathcal Z}_{\vec{d}}^{\Gamma}(q)
=(-1)^{\sum_{i=1}^r \gamma_id_i}
\sum_{(\mu^i,\alpha^i,\beta^i)}
\frac{1}{z_{\vec{\mu}}z_{\vec{\alpha}}z_{\vec{\beta}}
[\vec{\alpha}][\vec{\beta}]}
\prod_{i=1}^r
\langle
\mu^{i-1}\cup\alpha^i|q^{-(\gamma_i+2){\mathcal F}_2}
      |\mu^{i}\cup\beta^i
\rangle.
\end{equation}
(For any tuple of partitions 
$\vec{\lambda}=(\lambda^1,\ldots,\lambda^l)$,
$z_{\vec{\lambda}}=z_{\lambda^1}\ldots z_{\lambda^l}$
and
$[\lambda]=[\lambda^1]\ldots[\lambda^l]$.)

(iii)
Let $\Gamma$ be the GT graph in
example \ref{expart}-3.
We take the following flag order $\iota$:
\begin{equation}\notag
\iota(v_1,e_1)=\iota(v_2,e_2)=\iota(v_3,e_3)=1,\qquad
\iota(v_4,e_4)=2.
\end{equation}
%

We write 
$\vec{\mu}=(\mu^1,\mu^2,\mu^3,\mu^4)$ 
where $\mu^i$ is a partition associated to $v_i$.
We also write
$\Vec{\nu}=(\nu^{1,1},\nu^{1,2},\nu^{1,4},\nu^{2,2},\nu^{2,3},
\nu^{3,3},\nu^{3,1},\nu^{4,4})$
where $\nu^{i,j}$ is a partition associated to the flag
$(v_i,e_j)$.
The conditions for
a {\em $\Gamma$-set} 
$(\Vec{\mu},\vec{\nu})$ to be of
degree $\Vec{d}=(d_1,d_2,d_3,d_4)$
are
\begin{equation}\notag
\begin{split}
&|\mu^1|+|\nu^{1,1}|=|\mu^3|+|\nu^{3,1}|=d_1,
\qquad
|\mu^2|+|\nu^{2,2}|
=|\mu^1|+|\nu^{1,2}|=d_2,
\\
&
|\mu^3|+|\nu^{3,3}|
=|\mu^2|+|\nu^{2,3}|=d_3,
%
\qquad|\nu^{1,4}|=|\nu^{4,4}|=d_4,
\qquad |\mu^4|=0.
\end{split}
\end{equation}
The partition function ${\mathcal Z}_{\Vec{d}}^{\Gamma}(q)$
is written as:
\begin{equation}\label{eq:zop}
\begin{split}
{\mathcal Z}_{\Vec{d}}^{\Gamma}(q)
&=
(-1)^{\sum_{i=1}^3(b_i+1)d_i+b_4d_4}
\sum_{\begin{subarray}{c}
(\Vec{\mu},\Vec{\nu});\text{$\Gamma$-set}\\
\text{of degree $\Vec{d}$}
\end{subarray}
}
\frac{(-1)^{l(\nu^{1,2})+l(\nu^{1,4})}}
{z_{\Vec{\mu}}z_{\Vec{\nu}}
[\nu^{1,4}][\nu^{2,2}][\nu^{2,3}][\nu^{3,1}][\nu^{3,3}][\nu^{4,4}]}
\\
&\times
\langle \mu^1\cup\nu^{1,1}|
q^{-(b_1+1){\mathcal F}_2}|\mu^3\cup \nu^{3,1}\rangle
\langle \mu^2\cup\nu^{2,2}|
q^{-(b_2+1){\mathcal F}_2}|\mu^1\cup \nu^{1,2}\rangle
\\&\times
\langle \mu^3\cup\nu^{3,3}|
q^{-(b_3+1){\mathcal F}_2}|\mu^1\cup \nu^{2,3}\rangle
\langle
\nu^{1,4}|{\mathcal E}_0(-\nu^{1,1}){\mathcal E}_0(\nu^{1,2})
q^{-b_4{\mathcal F}_2}
|\nu^{4,4}\rangle.
\end{split}
\end{equation}
\end{example}

Now we describe a proof of lemma \ref{zop} in
the last example.

We apply lemma \ref{threepoint} to $C_{\Vec{\lambda}_{v_i}}(q)$
$(1\leq i\leq 4)$:
\begin{equation}\notag
\begin{split}
C_{\Vec{\lambda}_{v_1}}(q)
&=
\sum_{\mu^1}
\sum_{\nu^{1,1},\nu^{1,4},\nu^{1,2}}
\frac{(-1)^{l(\nu^{1,1})}}
{z_{\mu^1}z_{\nu^{1,1}}z_{\nu^{1,4}}z_{\nu^{1,2}}[\nu^{1,4}]}
\\&\times
\langle v_{\lambda^{1t}}|\nu^{1,1}\cup \mu^1\rangle
\langle v_{\lambda^{2t}}|q^{-{\mathcal F}_2}|\nu^{1,2}\cup \mu^1\rangle
\langle v_{\lambda^4}|{\mathcal E}_0(-\nu^{1,1})
                      {\mathcal E}_0(\nu^{1,2})|\nu^{1,4}\rangle.
\\
%
%
C_{\Vec{\lambda}_{v_2}}(q)&=
\sum_{\mu^2}\sum_{\nu^{2,2}\nu^{2,3}}
\frac{(-1)^{l(\nu^{2,2})}}{z_{\mu^2}z_{\nu^{2,2}}z_{\nu^{2,3}}}
\langle v_{\lambda^{2t}}|\nu^{2,2}\cup\mu^2 \rangle
\langle v_{\lambda^{3t}}|q^{-{\mathcal F}_2}|\nu^{2,3}\cup\mu^2\rangle
\langle 0|{\mathcal E}_0(-\nu^{2,2}){\mathcal E}_0(\nu^{2,3})|0\rangle
\\
&=\sum_{\mu^2}\sum_{\nu^{2,2}\nu^{2,3}}
\frac{1}{z_{\mu^2}z_{\nu^{2,2}}z_{\nu^{2,3}}[\nu^{2,2}][\nu^{2,3}]}
\langle v_{\lambda^{2t}}|\nu^{2,2}\cup\mu^2 \rangle
\langle v_{\lambda^{3t}}|q^{-{\mathcal F}_2}|\nu^{2,3}\cup\mu^2\rangle.
\\
%
%
C_{\Vec{\lambda}_{v_3}}(q)&=
\sum_{\mu^3}\sum_{\nu^{3,3}\nu^{3,1}}
\frac{1}{z_{\mu^3}z_{\nu^{3,3}}z_{\nu^{3,1}}[\nu^{3,3}][\nu^{3,1}]}
\langle v_{\lambda^{3t}}|\nu^{3,3}\cup\mu^3 \rangle
\langle v_{\lambda^{1t}}|q^{-{\mathcal F}_2}|\nu^{3,1}\cup\mu^3\rangle.
\\
%
C_{\Vec{\lambda}_{v_4}}(q)&=
\sum_{\nu^{4,4}}
\frac{1}{z_{\nu^{4,4}}[\nu^{4,4}]}
\langle v_{\lambda^{4t}}|\nu^{4,4}\rangle.
\end{split}
\end{equation}
The factor $q^{b_i\frac{\kappa(\lambda^i)}{2}}$ $(1\leq i\leq 4)$
is equal to
\begin{equation}\notag
q^{b_i\frac{\kappa(\lambda^i)}{2}}
=\langle
v_{\lambda^i}|q^{b_i{\mathcal F}_2}|v_{\lambda^i}
\rangle
=
\langle
v_{\lambda^{it}}|q^{-b_i{\mathcal F}_2}|v_{\lambda^{it}}
\rangle.
\end{equation}
Next we perform the summation over $\lambda^2$:
\begin{equation}\notag
\begin{split}
&\sum_{\lambda^2}
\langle 
v_{\lambda^{2t}}|q^{-b_2{\mathcal F}_2}|v_{\lambda^{2t}}\rangle
\langle v_{\lambda^{2t}}|q^{-{\mathcal F}_2}|\nu^{1,2}\cup \mu^1\rangle
\langle v_{\lambda^{2t}}|\nu^{2,2}\cup\mu^2 \rangle
\\
&=\langle \mu^2\cup \nu^{2,2}|q^{-(b_2+1){\mathcal F}_2}
                           |\nu^{1,2}\cup \mu^1 \rangle
\end{split}
\end{equation}
The summations over $\lambda^1$ and $\lambda^3$ are similar.
The summation over $\lambda^4$ is:
\begin{equation}\notag
\begin{split}
&\sum_{\lambda^4}
\langle v_{\lambda^4}|q^{b_4{\mathcal F}_2}|v_{\lambda^4} \rangle
\langle v_{\lambda^4}|{\mathcal E}_0(-\nu^{1,1})
                      {\mathcal E}_0(\nu^{1,2})|\nu^{1,4}\rangle
\langle v_{\lambda^{4t}}|\nu^{4,4}\rangle
\\
&=
\langle \nu^{1,4}|{\mathcal E}_0(-\nu^{1,1})
                      {\mathcal E}_0(\nu^{1,2})
q^{b_4{\mathcal F}_2} 
R\,
|\nu^{4,4}\rangle
\\
&=
(-1)^{d_4+l(\nu^{1,4})+l(\nu^{1,1})+l(\nu^{1,2})}
\langle \nu^{1,4}|{\mathcal E}_0(-\nu^{1,2})
                      {\mathcal E}_0(\nu^{1,1})
q^{-b_4{\mathcal F}_2} 
|\nu^{4,4}\rangle.
\end{split}
\end{equation}
In the middle line, $R$ is the transposition operator
introduced in subsection \ref{trop}.
In passing to the last line,
we have moved $R$ to the left using
(\ref{ar}) (\ref{aar})
and then exchanged ${\mathcal E}_0(\nu^{1,1})$
and ${\mathcal E}_0(\nu^{1,2})$
using the fact that ${\mathcal E}_0(i)$ and 
${\mathcal E}_0(j)$ $(i,j\in\mathbb{Z}\setminus \{0\})$
commute with each other.
Combining the above expressions, 
we obtain (\ref{eq:zop}).

Lemma \ref{zop} is proved for a general GT graph
completely in the same manner.

\subsection{Graph Expression}\label{ge}
In this subsection, 
we will express 
${\mathcal Z}_{\vec{d}}^{\Gamma}(q)$
as the sum over a certain set of labeled graphs.

Before proceeding,
we explain briefly the graph expression 
introduced in \cite{ko},
section 3.2.
\subsubsection{Graph Expression of VEV}
Let $\Vec{c}$ and $\vec{n}$ be 
sequences of integers of the same length $l$:
\begin{equation}\notag
\Vec{c}=(c_1,\ldots,c_l),\qquad
\Vec{n}=(n_1,\ldots,n_l)
\end{equation}
such that
$(c_i,n_i)\neq (0,0)$ for $1\leq i\leq l$)
and
$|\Vec{c}|:=\sum_{i=1}^lc_i=0$.

The  vacuum expectation value
\begin{equation}\label{vev}
\langle {\mathcal E}_{c_1}(n_1)
{\mathcal E}_{c_2}(n_2)\cdots
{\mathcal E}_{c_l}(n_l)
\rangle
\end{equation}
is computed by applying the commutation relation:
\begin{equation}\notag
{\mathcal E}_a(m){\mathcal E}_b(n)=
{\mathcal E}_b(n){\mathcal E}_a(m)
+
\begin{cases}
[an-bm]\,
{\mathcal E}_{a+b}(m+n)
& (a+b,m+n)\neq (0,0)\\
a&(a+b,m+n)=(0,0)
\end{cases}
\end{equation}
For the algorithm to be well-defined,
we set the rule that the commutation relation is applied to the 
rightmost neighboring 
pair $({\mathcal E}_a(m),{\mathcal E}_b(n))$
such that $a\geq 0$ and $b<0$.
We also use
the relations:
\begin{equation}\notag
\begin{split}
&\langle\cdots{\mathcal E}_a(m)\rangle
=
\begin{cases}0&(a>0)\\ 
\langle\cdots\rangle \frac{1}{[m]}&(a=0)
\end{cases}
\\
&\langle{\mathcal E}_b(n)\cdots\rangle=0\quad (b<0),
\qquad \text{ and }\langle 1\rangle=1.
\end{split}
\end{equation}
We associate to the commutation relation the drawing:
%
\newcommand{\before}{  
\unitlength .15cm
\begin{picture}(7,8)(-3,-2)
\thicklines
\put(0,5){\circle*{1}}
\put(5,5){\circle*{1}}
\put(-3,7){$(a,m)$}
\put(3,7){$(b,n)$}
\end{picture}
}
\newcommand{\crossing}{
\unitlength .15cm
\begin{picture}(7,8)(-3,-2)
\thicklines
\put(0,5){\line(1,-1){5}}
\put(5,5){\line(-1,-1){5}}
\put(0,5){\circle*{1}}
\put(5,5){\circle*{1}}
\put(-3,7){$(a,m)$}
\put(3,7){$(b,n)$}
\put(0,0){\circle*{1}}
\put(5,0){\circle*{1}}
\put(-3,-2){$(b,n)$}
\put(3,-2){$(a,m)$}
\end{picture}
}
\newcommand{\collide}{
\unitlength .15cm
\begin{picture}(7,8)(-3,-3)
\thicklines
\put(0,5){\line(1,-2){2.5}}
\put(5,5){\line(-1,-2){2.5}}
\put(0,5){\circle*{1}}
\put(5,5){\circle*{1}}
\put(-3,7){$(a,m)$}
\put(3,7){$(b,n)$}
\put(2.5,0){\circle*{1}}
\put(-2,-2){$(a+b,m+n)\neq (0,0)$}
\end{picture}
}
\newcommand{\coll}{
\unitlength .15cm
\begin{picture}(7,10)(-3,-2)
\thicklines
\put(0,5){\line(1,-2){2.5}}
\put(5,5){\line(-1,-2){2.5}}
\put(0,5){\circle*{1}}
\put(5,5){\circle*{1}}
\put(-3,7){$(a,m)$}
\put(3,7){$(b,n)$}
\put(2.5,0){\circle{1}}
\put(-1,-2){$(a+b,m+n)=(0,0)$}
\end{picture}
}
%
\begin{equation}\notag
\raisebox{-.7cm}{\before} 
\qquad = \qquad
\raisebox{-.7cm}{\crossing} 
\qquad + \qquad
\begin{cases}
\raisebox{-.7cm}{\collide}\\
\\
\raisebox{-.5cm}{\coll}
\end{cases}
\end{equation}
%
%
Then graphs are generated 
over the course of the calculation.
%
\begin{definition}%
${\rm Graph}^{\bullet}(\vec{c},\vec{n})$
is
the set of graphs generated by this procedure.
%
\end{definition}%

Every graph $F\in{\rm Graph}_a^{\bullet}(\vec{c},\vec{n}) $ is 
the graph union of
(binary rooted) trees.
We call $F$ a {\em VEV forest}
and 
a connected component of $F$ a {\em VEV tree}.
The two component label of
every vertex $v$ 
is denoted by $(c_v,n_v)$ and called the 
{\em vertex label} of $v$.
We add another label called the {\em leaf index}
to every leaf of $F$:
leaves of $F$ correspond one-to-one  to 
components $(c_i,n_i)$ of $(\vec{c},\vec{n})$;
so the leaf index of
a leaf corresponding to $(c_i,n_i)$
is defined to be $i$.

By construction, every graph 
$F\in {\rm Graph}^{\bullet}(\vec{c},\vec{n})$
represents one term in the final result of the
computation of the VEV (\ref{vev}).
The corresponding term can be recovered from $F$ as follows.
For each VEV tree $T$ in $F$,
let $V_2(T)$ be the set of vertices
which have two adjacent vertices
at the upper level.
The upper left and right
vertices adjacent to $v\in V_2(T)$
are denoted by
$L(v)$ and $R(v)$, respectively.
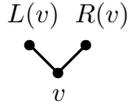
\newcommand{\LR}{
\unitlength .15cm
\begin{picture}(10,8)(0,0)
\thicklines
\put(5,5){\line(1,-1){2.5}}
\put(10,5){\line(-1,-1){2.5}}
\put(3,7){$L(v)$}
\put(5,5){\circle*{1}}
\put(10,5){\circle*{1}}
\put(9,7){$R(v)$}

\put(7.5,2.5){\circle*{1}}
\put(7,0){$v$}
\end{picture}
}
\begin{equation}\notag
\raisebox{-.5cm}{\LR}
\end{equation}
For a vertex $v\in V_2(T)$, 
define
\begin{equation}\notag
\xi_v=
c_{L(v)}n_{R(v)}-c_{R(v)}n_{L(v)}.
\end{equation}
We define the {\em amplitude} ${\mathcal A}(F)$ by
\begin{equation}\notag
\begin{split}
{\mathcal A}(F)&=
\prod_{T:\text{ VEV tree in $F$}}
{\mathcal A}(T),
\\
{\mathcal A}(T)&=
\begin{cases}
{\prod_{v\in V_2(T)}[\xi_v]}/
{[n_{\text{root}}]}
&(n_{\text{root}}\neq 0)
\\
c_{L(\text{root})}
\prod_{\begin{subarray}{c}v\in V_2(T),\\
     v\neq \text{root}\end{subarray}}
[\xi_v]
&(n_{\text{root}}=0)
\end{cases}
\end{split}
\end{equation}
%
%
The amplitude ${\mathcal A}(F)$
is exactly
the term corresponding to $F$.
Thus we have
\begin{prop}\label{prop:vev}
\begin{equation}\notag
\langle {\mathcal E}_{c_1}(n_1)
{\mathcal E}_{c_2}(n_2)\cdots
{\mathcal E}_{c_l}(n_l)
\rangle
=\sum_{F\in {\rm Graph}^{\bullet}(\Vec{c},\Vec{ n})}
{\mathcal A}(F).
\end{equation}
\end{prop}
\begin{proof} Clear.
\end{proof}

The amplitude of a VEV tree $T$
admits an important pole structure.
Let us define $m_v:=\textrm{gcd}(c_v,n_v)$ for a leaf $v$
of a VEV tree $T$
and 
\begin{equation}\label{def:B}
m(T):=\textrm{gcd}(\{m_v\}_{v:\text{leaf of $T$}}),
\qquad
{\mathcal B}(T):=\frac{{\mathcal A}(T)}{\prod_{v:\text{leaf}}[m_v]}.
\end{equation}
%
\begin{prop}\label{poleB}
Let $T$ be a VEV tree.
\begin{enumerate}
\item %
There exists $g_{T}\in \mathbb{Z}$ and
$f_{T}(t)\in \mathbb{Z}[t]$ such that
\begin{equation}\notag
{\mathcal B}(T)=\begin{cases}
\frac{g_{T}}{t_{m(T)}}+f_{T}(t_{m(T)})& 
(\text{$m(T)$ odd or $n_{\text{root}}/m(T)$ even})
\\
\frac{g_{T}}{t_{m(T)}}
\Big(1+\frac{t_{m(T)/2}}{2}\Big)+f_{T}(t_{m(T)})& 
(\text{$m(T)$ even and $n_{\text{root}}/m(T)$ odd})
\end{cases}
\end{equation}
\item %
Moreover,
\begin{equation}\notag
g_{T}=g_{T_{(0)}}\cdot {m(T)}^{\#\text{leaves}-1}.
\end{equation}
Here  $T_{(0)}$ is the VEV tree 
which is the 
same as $T$
except all the vertex-labels are multiplied by $1/m(T)$.
\end{enumerate}
\end{prop}%
See proposition 6.1 and section 6.3.2 \cite{ko} for a proof.
%
\subsubsection{Graph Expression}
Now we will express 
${\mathcal Z}_{\vec{d}}^{\Gamma}(q)$
as the sum over a certain set of labeled graphs.
This takes two steps.

The first step is to write 
the matrix elements $K_e(\vec{\mu},\vec{\nu})$ 
as the sum over a set of labeled graphs using 
proposition \ref{prop:vev}.

For partitions $\eta,\xi,\xi',\eta'$ such that 
$|\eta|=|\eta'|$,
we define
${\rm Graph}_a^{\bullet}(\eta,-\xi,\xi',\eta')$
to be the set
${\rm Graph}^{\bullet}(\Vec{c},\Vec{n})
$
with 
\begin{equation}\notag
\begin{split}
\Vec{c}&=(\eta_{l(\eta)},\ldots,\eta_1,
\underbrace{0,\ldots,0}_{l(\xi)+l(\xi')},
-\eta_1',\ldots,-\eta'_{l(\eta')}),
\\
\Vec{n}&=(\underbrace{0,\ldots,0}_{l(\eta)},
-\xi_1,\ldots,-\xi_{l(\xi)},
\xi_1',\ldots,\xi'_{l(\xi')},
a\eta_1',\ldots,a\eta'_{l(\eta')}).
\end{split}
\end{equation}
%
\begin{definition}
We set
\begin{equation}\notag
{\rm Graph}_e^{\bullet}(\Vec{\mu},\Vec{\nu})
:={\rm Graph}_{N_e}^{\bullet}(\eta_e,-\xi_e,\xi_e',\eta_e')
\qquad (e\in E_3(\Gamma)).
\end{equation}
Here $\eta_e,\xi_e,\xi_e',\eta_e'$ are partitions
listed in table \ref{tab:xi}.
%
\end{definition}

\begin{table}[t]
\begin{equation}\notag
\begin{array}{|c|cccc|}\hline
(\iota(f),\iota(f'))& \eta_e & \xi_e & \xi_e' &\eta_e' 
\\\hline
(1,1)&\mu\cup \nu& \emptyset&\emptyset& \mu'\cup\nu'
\\
(1,2)&\mu\cup\nu&\nu^{\prime 1}&\nu^{\prime 3}&\nu'
\\
(1,3)&\mu\cup \nu&\emptyset&\emptyset&\mu'\cup\nu'
\\
(2,1)&\nu &\nu^3&\nu^1&\mu'\cup\nu'
\\
(2,2)&\nu&\nu^3\cup \nu^{\prime 1}&\nu^{1}\cup\nu^{\prime 3}&
\nu'
\\
(2,3)&\nu&\nu^1&\nu^3&\mu'\cup\nu'
\\
(3,1)&\mu\cup\nu&\emptyset&\emptyset&\mu'\cup\nu'
\\
(3,2)&\mu\cup\nu&\nu^{\prime 1}&\nu^{\prime 3}&\nu'
\\
(3,3)&\mu\cup\nu&\emptyset&\emptyset&\mu'\cup\nu' 
\\\hline
\end{array}
\end{equation}
\caption{$\eta_e,\xi_e,\xi_e',\eta_e'$ for $e\in E_3(\Gamma)$.
The notation is the same as in table  \ref{tab:NL}.}
\label{tab:xi}
\end{table}
%

Then the matrix element $K_e(\vec{\mu},\vec{\nu})$
is written as follows:
\begin{lemma}\label{lem:Ke} 
\begin{equation}\notag
K_e(\vec{\mu},\vec{\nu})=
\sum_{F\in {\rm Graph}_e^{\bullet}(\vec{\mu},\vec{\nu})}
{\mathcal A}(F).
\end{equation}
\end{lemma}
\begin{proof}
Since 
(\cite{ko}, eq (3))
\begin{equation}\notag
q^{a{\mathcal F}_2}\alpha_{-i}q^{-a{\mathcal F}_2}
={\mathcal E}_{-i}(ai),
\end{equation}
a matrix element of the form
\begin{equation}\notag
\langle \eta |{\mathcal E}_0(-\xi){\mathcal E}_0(\xi') 
q^{a{\mathcal F}_2}|\eta' \rangle
\qquad (\eta,\xi,\xi',\eta':\text{partitions}, 
|\eta|=|\eta'|)
\end{equation}
is equal to the VEV 
\begin{equation}\notag
\langle
{\mathcal E}_{\eta_{l(\eta)}}(0)\ldots
{\mathcal E}_{\eta_{1}}(0)
\,
{\mathcal E}_0(\xi)
\,
{\mathcal E}_0(-\xi)
\,
{\mathcal E}_{-\eta'_1}(a\eta'_1)
\ldots
{\mathcal E}_{-\eta'_{l(\eta')}}(a\eta'_{l(\eta')})
\rangle.
\end{equation}
Therefore the lemma follows from proposition \ref{prop:vev}.
\end{proof}

The second step is to construct a new type of 
labeled graphs  
so that
${\mathcal Z}_{\vec{d}}^{\Gamma}(q)$
becomes the sum of quantities over these graphs.

Let us rewrite ${\mathcal Z}_{\vec{d}}^{\Gamma}(q)$
using ${\mathcal B}(F)$ introduced in
(\ref{def:B}):
$
{\mathcal B}(F)={\mathcal A}(F)/[\eta][\xi][\xi'][\eta']
$
for $F\in{\rm Graph}^{\bullet}_a(\eta,-\xi,\xi',\eta')$.
By lemmas \ref{zop} and \ref{lem:Ke},
\begin{equation}\notag
{\mathcal Z}_{\vec{d}}^{\Gamma}(q)
=
\sum_{\begin{subarray}{c}
  (\vec{\mu},\vec{\nu});\text{ $\Gamma$-set}\\
   \text{of degree $\vec{d}$}\end{subarray}}
\frac{(-1)^{L_1(\vec{\mu},\vec{\nu})+L_2(\vec{\mu},\vec{\nu})}
}{z_{\vec{\mu}}z_{\vec{\nu}}}
\underbrace{
\prod_{v\in V_3(\Gamma)}
[\mu^v]^2
\prod_{
         f\in \hat{F}_3(\Gamma)}
[\nu^f]^2
}_{\ast}
\sum_{(F_e)}
\prod_{e\in E_3(\Gamma)}
{\mathcal B}(F_e).
\end{equation}
Here 
\begin{equation}\notag        
\hat{F}_3(\Gamma):=\{f\in F_3(\Gamma)|
\iota(f)=1,3,
         \myf_2(f)\in F_3(\Gamma)\},
\end{equation}
and
the second summation is over the set
\begin{equation}\notag
\prod_{e\in E_3(\Gamma)}
{\rm Graph}_e^{\bullet}(\vec{\mu},\vec{\nu}).
\end{equation}
So we need to incorporate the  factors in $(\ast)$.
Recall that, by construction,
leaves in
a VEV forest $F_e\in
{\rm Graph}_e^{\bullet}(\vec{\mu},\vec{\nu})$ 
$(e\in E_3(\Gamma))$
correspond one-to-one
to parts of
the partitions $\eta_e,\xi_e,\xi_e',\eta_e'$.
Therefore,
every leaf corresponds to a part of $\mu^v$ $(v\in V_3(\Gamma))$
or $\nu^f$ $(f\in F_3(\Gamma))$ ${}^\dagger$
(cf. remark \ref{leaves}).
For $\mu^v$, 
there are two leaves
associated to every part $\mu^v_i$,
the one in $F_e$ such that $(v,e)=\myf_1(v)$
and the other in $F_{e'}$ such that $(v,e')=\myf_3(v)$.
Similarly,
for $\nu^f$ with $f\in \hat{F}_3(\Gamma)$,
there are two leaves associated to 
every part $\nu^f_i$ of $\nu^f$ with
$f\in \hat{F}_3(\Gamma)$,
the one in  $F_e$ such that $f=(v,e)$
and the other in $F_{e'}$ such that
$\myf_2(f)=(v,e')$.

We construct a new graph from
\begin{equation}\notag
(F_e) 
\in 
\prod_{e\in E_3(\Gamma)}{\rm Graph}^{\bullet}_{e}(\Vec{\mu},\Vec{\nu})
\end{equation}
as follows.
\begin{enumerate}
\item
Assign the label $e\in E_3(\Gamma)$ 
to each $F_e$ and make the graph union.
\item
Join the
two leaves associated to 
$\mu^v_i$ $(v\in V_3(\Gamma), 1\leq i\leq l(\mu^v))$
and attach the label $\mu^v_i$ to the new edge.
Also join the 
two leaves associated to
$\nu^f_i$ $(f\in \hat{F}_3(\Gamma),1\leq i\leq l(\nu^f))$
and attach the label $\nu^f_i$ to the new edge.
\end{enumerate}

The resulting graph $W$ is a set of VEV forests marked by 
$e\in E_3(\Gamma)$
and
joined through leaves.
We call $W$ a {\em combined forest}.
The new edges are called the {\em bridges}.
The label of a bridge $b$ is denoted by $h(b)$.
%
\begin{definition}%
The set of combined forests constructed by the above procedure
is denoted by
${\rm Comb}^{\bullet}_{\Gamma}
(\Vec{\mu},\Vec{\nu})$.
The subset consisting of connected combined forests
is denoted by
${\rm Comb}_{\Gamma}^{\circ}
(\Vec{\mu},\Vec{\nu})$.
\end{definition}%
Examples of combined forests are shown in figure \ref{excomb}.

\begin{figure}[tt]
\begin{equation}\notag
\hspace*{-5cm}
\ancomb{\treea}{\treea}
\end{equation}
\begin{equation}\notag
\hspace*{-5cm}
\ancomb{\treec}{\treeb}
\end{equation}
$\Gamma$: the GT graph in example \ref{expart}-3
\\
$b_1,b_2,b_3\neq -1$ and $b_4\neq 0$.
\\
$\mu^1=\mu^2=\mu^3=(1)$, $\mu^4=\emptyset$,
\\
$\nu^{1,1}=\nu^{1,4}=\nu^{2,3}
=\nu^{3,3}=\nu^{3,1}=\nu^{4,4}=(1)$,
$\nu^{1,2}=\nu^{2,2}=\emptyset$.
\caption{
Example of combined forests.
(Vertex labels and leaf indices are omitted.)
}
\label{excomb}
\end{figure}

For a combined forest
$W
\in {\rm Comb}_{\Gamma}^{\bullet}
(\vec{\mu},\vec{\nu})$,
we define
\begin{equation}\notag
{\mathcal H}(W)=
(-1)^{L_1(W)+L_2(W)}
\prod_{e\in E_3(\Gamma)}{\mathcal B}(F_e)
\prod_{b;\text{bridge}}[h(b)]^2,
\end{equation}
where 
$L_1(W)=L_1(\vec{\mu},\vec{\nu})$ and
$L_2(W)=L_2(\vec{\mu},\vec{\nu})$.
${\mathcal H}(W)$ is called the {\em combined amplitude}.

\begin{prop}%
%
\begin{equation}\notag
{\mathcal Z}^{\Gamma}_{\Vec{d}}(q)
=\sum_{\begin{subarray}{c}
(\Vec{\mu},\Vec{\nu});\\
\text{$\Gamma$-set of}\\
\text{degree $\Vec{d}$}
\end{subarray}}
\frac{1}
{z_{\Vec{\mu}}z_{\Vec{\nu}}}
\sum_{W\in {\rm Comb}^{\bullet}_{\Gamma}
(\Vec{\mu},\Vec{\nu})}
{\mathcal H}(W).
\end{equation}
\end{prop}
\begin{proof}
The proposition follows from the definitions of 
the combined forest and the combined amplitude.
\end{proof}

\begin{remark}\label{leaves}
Precisely speaking, 
the statement ($\dagger$) is not correct.
As an example, consider 
a VEV forest $F_e\in {\rm Graph}_e^{\bullet}
(\vec{\mu},\vec{\nu}) $
where $e\in E_3(\Gamma)$ such that
$\eta_e=\mu^v\cup\nu^f$ with $f=(v,e)$.
If
$\mu^v$ and $\nu^f$ have equal parts,
say $\mu^v_i=\cdots=\mu^v_{i+l-1}=\nu^f_j=\cdots=\nu^f_{j+m-1}
=k\in \mathbb{N}$,
then 
$F_e$
has
$l+m$ leaves with the vertex label $(k,0)$.
But there is no preferred way to
determine which leaves correspond to which parts.
To solve
this problem, 
we promise that
these leaves correspond, from left to right,
to 
$\mu^v_{i+l-1},\ldots,\mu^v_i,\nu^f_{j+m-1},\ldots,\nu^f_j$
(ie, the outer $l$ leaves correspond to the parts of $\mu^v$,
and the inner $m$ leaves to those of $\nu^f$).
Similarly,
if $\eta_e'=\mu^v\cup\nu^f$ and 
$\mu^v$ and $\nu^f$ have equal parts,
we promise that outer leaves correspond to the parts
of $\mu^v$ and inner leaves to the parts of $\nu^f$.
If $\xi_e=\nu^{\myf_3(v)}\cup \nu^{\myf_1(v')}$
(resp. $\xi_e'=\nu^{\myf_1(v)}\cup \nu^{\myf_3(v')}$)
and if $\nu^{\myf_3(v)}$ and $\nu^{\myf_1(v')}$
(resp. $\nu^{\myf_1(v)}$ and $\nu^{\myf_3(v')}$)
 have equal parts,
we promise that left leaves correspond to the parts of
$\nu^{\myf_3(v)}$ (resp. $\nu^{\myf_1(v)}$)
and right leaves to the parts of 
$\nu^{\myf_1(v')}$ (resp. $\nu^{\myf_3(v')}$).
\end{remark}

\subsection{Free Energy}\label{feen}
We take the 
logarithm of the partition function by using
the exponential formula.
As a result, we obtain
the free energy as
the sum over
connected combined graphs.
\begin{prop}\label{free}
%
\begin{equation}\notag
{\mathcal F}^{\Gamma}_{\Vec{d}}(q)
=
\sum_{\begin{subarray}{c}
(\Vec{\mu},\Vec{\nu});\\
\text{$\Gamma$-set of}\\
\text{degree $\Vec{d}$}
\end{subarray}}
\frac{1}
{z_{\Vec{\mu}}z_{\Vec{\nu}}}
\sum_{W\in {\rm Comb}^{\circ}_{\Gamma}
(\Vec{\mu},\Vec{\nu})}
{\mathcal H}(W).
\end{equation}
\end{prop}%
\begin{proof}
We use the formulation in \cite{ko},
appendix A.

For a combined forest $W$, we define $\overline{W}$ to be 
the graph obtained by forgetting all leaf indices.
Two  combined forests $W$ and $W'$ are {\em equivalent}
if $\overline{W}$ and $\overline{W}'$ are isomorphic
as labeled graphs.
The set of equivalence classes in
${\rm Comb}^{\bullet}_{\Gamma}(\Vec{\mu},\Vec{\nu})$ 
is  denoted by
$\overline{{\rm Comb}}^{\bullet}_{\Gamma}
(\Vec{\mu},\vec{\nu})$.
We define
\begin{equation}\notag
\begin{split}
\overline{{\rm Comb}}^{\bullet}_{\Gamma}(d)&=
\coprod_{|\Vec{d}|=d}
\coprod_{
\begin{subarray}{c}
(\Vec{\mu},\Vec{\nu})\\
\text{$\Gamma$-set of}\\
\text{degree $\Vec{d}$}
\end{subarray}}
\overline{{\rm Comb}}^{\bullet}_{\Gamma}
(\Vec{\mu},\Vec{\nu}) \qquad (d\geq 1),
\\
\overline{{\rm Comb}}^{\bullet}_{\Gamma}
&=\coprod_{d\geq 1}
\overline{{\rm Comb}}_{\Gamma}^{\bullet}(d).
\end{split}
\end{equation}
The set $\overline{{\rm Comb}}^{\bullet}_{\Gamma}$ is a G-set.

Next we define a map 
$\Psi:{\rm Comb}_{\Gamma}^{\bullet}\to \mathbb{Q}(t)[[\vec{Q}]]$,
so that it satisfies
\begin{equation}\label{def:Psi}
\frac{1}{\#\aut(G)}\Psi(\overline{G})
=
\sum_{\begin{subarray}{c}
W\in{\rm Comb}_{\Gamma}^{\bullet}(\vec{\mu},\vec{\nu})\\
\overline{W}=G
\end{subarray}
}
\frac{1}{z_{\Vec{\mu}}z_{\Vec{\nu}}}\,
{\mathcal H}(W) \,
\Vec{Q}^{\Vec{d}}
\qquad
(G\in 
\overline{{\rm Comb}}_{\Gamma}^{\bullet}(\Vec{\mu},\Vec{\nu})
)
\end{equation}
where
$\Vec{d}$ is the degree of the $\Gamma$-set
$(\vec{\mu},\vec{\nu})$.
It is easy to see that
$\Psi$ 
is a grade preserving map.
Showing that
$\Psi$ is multiplicative with respect to the graph union
is straightforward.
See appendix \ref{appendix:Psi}.

$(\overline{{\rm Comb}}^{\bullet}_{\Gamma},
\mathbb{Q}(t)[[\vec{Q}]],\Psi)$
satisfies the conditions for the GA-triple.
Applying
the exponential formula to this,
we obtain
\begin{equation}\notag
\log\Big[1+\sum_{G\in   \overline{{\rm Comb}}_{\Gamma}^{\bullet}}
\frac{1}{\#\aut(G)}
\Psi(G)\Big]
=
\sum_{G\in   \overline{{\rm Comb}}_{\Gamma}^{\circ}}
\frac{1}{\#\aut(G)}
\Psi(G).
\end{equation}
Here
 $\overline{{\rm Comb}}_{\Gamma}^{\circ}$
stands for the subset of
$\overline{{\rm Comb}}_{\Gamma}^{\bullet}$
consisting of all connected graphs.
Since the LHS is equal to
$\log {\mathcal Z}^{\Gamma}(q,\Vec{Q})$,
the RHS is equal to
the free energy
${\mathcal F}^{\Gamma}(q,\Vec{Q})$.
Substituting (\ref{def:Psi}) back, we obtain
proposition \ref{free}.
\end{proof}
\subsection{Proof of Proposition \ref{polynomiality}}
\label{ppp}
Finally we will give a proof of proposition \ref{polynomiality}.

Given a combined forest $W$ and a positive integer $k$,
there exists a unique
combined forest 
which is 
the same as $W$ 
except that all the vertex-labels are multiplied by $k$.
It is denoted by
$W_{(k)}$.

We first rewrite the free energy as follows.
\begin{equation}\notag
{\mathcal F}_{\Vec{d}}^{\Gamma}(q)=
\sum_{k;k|d_0}
\sum_{\begin{subarray}{c}
(\Vec{\mu},\Vec{\nu});\\
\text{$\Gamma$-set of}\\
\text{degree $\Vec{d}/k$};\\
{\rm gcd}(\Vec{\mu},\Vec{\nu})=1
\end{subarray}}
\frac{1}
{z_{k\Vec{\mu}}z_{k\Vec{\nu}}}
\sum_{W\in {\rm Comb}^{\circ}_{\Gamma}
(\Vec{\mu},\Vec{\nu})}
{\mathcal H}(W_{(k)})
\end{equation}
where $d_0={\rm gcd}(\Vec{d})$.
Then,
${G}_{\Vec{d}}^{\Gamma}(q)$ is equal to
\begin{equation}\notag
\begin{split}
{G}_{\Vec{d}}^{\Gamma}(q)&=
\sum_{k;k|{d}_0}
\sum_{\begin{subarray}{c}
(\Vec{\mu},\Vec{\nu});\\
\text{$\Gamma$-set of}\\
\text{degree $\Vec{d}/k$},\\
{\rm gcd}(\Vec{\mu},\Vec{\nu})=1
\end{subarray}}
\sum_{k';k'|k}
\frac{k'}{k}
\mu
\Big(\frac{k}{k'}\Big)
\frac{1}
{z_{k'\Vec{\mu}}z_{k'\Vec{\nu}}}
\sum_{W\in {\rm Comb}^{\circ}_{\Gamma}
(\Vec{\mu},\Vec{\nu})}
{\mathcal H}(W_{(k')})|_{q\to q^{k/k'}}
\\
&=\sum_{k;k|d_0}
\sum_{\begin{subarray}{c}
(\Vec{\mu},\Vec{\nu});\\
\text{$\Gamma$-set of}\\
\text{degree $\Vec{d}/k$};\\
{\rm gcd}(\Vec{\mu},\Vec{\nu})=1
\end{subarray}}
\frac{1}{k\#\aut(\Vec{\mu})\#\aut(\Vec{\nu})}
\sum_{W\in {\rm Comb}^{\circ}_{\gamma}(\Vec{\mu},\Vec{\nu})}
{\mathcal G}_k^{\Gamma}(W),
\end{split}
\end{equation}
where
\begin{equation}\notag
{\mathcal G}^{\Gamma}_k(W)=
\sum_{k':k'|k}
\mu\Big(\frac{k}{k'}\Big)
{{k'}^{-l(\Vec{\mu})-l(\Vec{\nu})+1}}
{\mathcal H}(W_{(k')})\Big|_{q\to q^{k/k'}}.
\end{equation}
Proposition \ref{polynomiality} follows from:
\begin{prop}
Let $(\vec{\mu},\vec{\nu})$ be a $\Gamma$-set 
such that
${\rm gcd}(\Vec{\mu},\Vec{\nu})=1$.
Let $W$ be a connected combined forest in 
${\rm Comb}^{\circ}_{\Gamma}(\Vec{\mu},\Vec{\nu})$
and let $k$ be a positive integer. Then
\begin{equation}\notag
t\cdot {\mathcal G}^{\Gamma}_k(W)\in\mathbb{Q}[t].
\end{equation}
\end{prop}
%
\begin{proof}
The proposition is a consequence of the next lemma.
\end{proof}
\begin{lemma}
Let $W$ be a connected combined forest in 
${\rm Comb}^{\circ}_{\Gamma}(\Vec{\mu},\Vec{\nu})$.
\begin{enumerate}
\item If $\#E(W)-\#V(W)+1>0$,
then ${\mathcal H}(W)\in \mathbb{Z}[t]$.
\item If $\#E(W)-\#V(W)+1=0$,
then $t_m{\mathcal H}(W)\in \mathbb{Z}[t]$
where $m=\textrm{gcd}(\vec{\mu},\vec{\nu})$.
Moreover, the followings hold:
\item
Assume that $\text{gcd}(\vec{\mu},\vec{\nu})=1$
and that $\#E(W)-\#V(W)+1=0$.
Then
\\
if $k$ is odd,
\begin{equation}\notag
{\mathcal H}(W_{(k)})-k^{l(\vec{\mu})+l(\vec{\nu})-1}
{\mathcal H}(W)|_{t\to t_k}\in \mathbb{Z}[t];
\end{equation}
if $k$ is even,
\begin{equation}\notag
{\mathcal H}(W_{(k)})-
k^{l(\vec{\mu})+l(\vec{\nu})-1}
{\mathcal H}(W)|_{t\to t_k}
\prod_{T\in VT_I(W)}\Big(1+\frac{t_{m(T)k/2}}{2}\Big)
\in \mathbb{Z}[t].
\end{equation}
Here $VT_I(W)$ is the set of VEV trees 
in $W$
such that $m(T)$ is odd and $n_{\text{root}}/m(T)$
is odd.
\end{enumerate}
\end{lemma}
\begin{proof}
This lemma follows from
proposition \ref{poleB}.
The proof is completely 
the same as those of propositions 6.7 and 6.8 in \cite{ko}.
\end{proof}%

This completes  the proof of proposition \ref{polynomiality}.
\appendix
\section{Proof of Lemma \ref{arith}}\label{proofarith}
In this appendix, we
give a proof of lemma \ref{arith}.
\subsection{Sublemmas}
Let $p$ be a prime integer, $i,b\in\mathbb{N}$.

\begin{lemma}\label{xpb}
\begin{equation}\notag
(x_1+\cdots+x_l)^{p^ib}
=
(x_1^p+\cdots+x_l^p)^{p^{i-1}b}+\text{ term divisible by }p^i
\qquad (i\geq 1)
\end{equation}
where $x_1,\ldots,x_l$ are
indeterminate variables. 
\end{lemma}
\begin{proof}
The proof is by induction on $i$.\\
$i=1$:
\begin{equation}\notag
\begin{split}
(x_1+\cdots +x_l)^{pb}&=
\big((x_1+\cdots+x_l)^p\big)^b
\\
&=(x_1^p+\cdots+x_l^p+
\text{ term divisible by }p)^b
\\
&=(x_1^p+\cdots+x_l^p)^b+
\text{ term divisible by }p).
\end{split}
\end{equation}
$i>1$: assume that this is true for $i-1$.
\begin{equation}\notag
\begin{split}
(x_1+\cdots+x_l)^{p^{i+1}b}
&=
\big((x_1+\cdots+x_l)^{p^ib}\big)^p
\\
&=
\big((x_1^p+\cdots +x_l^p)^{p^{i-1}b}+
\text{ term divisible by }p^i
\big)^p
\\
&=
\big
(x_1^p+\cdots +x_l^p)^{p^{i}b}
\text{ term divisible by }p^{i+1}.
\end{split}
\end{equation}
So the lemma is proved.
\end{proof}
\begin{lemma}\label{art}
Let
$n=(n_1,\ldots, n_l)$ be the set of nonnegative integers
such that
${\rm gcd}(n)=1$.
\begin{enumerate}
\item
For $k\in\mathbb{N}$,
\begin{equation}\notag
\frac{(k|n|)!}{(kn_1)!\cdots(kn_l)!}\equiv0 \mod |n|.
\end{equation}
\item
For a positive integer $k\in\mathbb{N}$ prime to $p$,
\begin{equation}\notag
\frac{(p^ik|n|)!}
{(p^ikn_1)!\cdots(p^ikn_l)!}
\equiv
\frac{(p^{i-1}k|n|)!}
{(p^{i-1}kn_1)!\cdots(p^{i-1}kn_l)!}
\mod p^i|n|.
\end{equation}
\end{enumerate}
\end{lemma}
\begin{proof}
1. 
The main idea of the proof is the following.
If a rational number $r$ satisfies
\begin{equation}\notag
h_1\cdot r\in\mathbb{Z}
\qquad\text{ and }
\qquad
h_2\cdot r\in\mathbb{Z}
\end{equation}
where $h_1,h_2$ are natural numbers,
then clearly it holds that
\begin{equation}\notag
{\rm gcd}(h_1,h_2)\cdot r\in\mathbb{Z}.
\end{equation}
Let us define
\begin{equation}\notag
C(k,n):=\frac{1}{|n|}
\frac{(k|n|)!}{(kn_1)!\ldots(kn_l)!}\in\mathbb{Z}.
\end{equation}
Then 
\begin{equation}\notag
n_1 C(k,n)=
\frac{(k|n|-1)!}{(kn_1-1)!\ldots(kn_l)!}\in\mathbb{Z}.
\end{equation}
The same holds for every $n_i$ ($1\leq i\leq l$).
Therefore
\begin{equation}\notag 
{\rm gcd}(n)\cdot  C(k,n)\in\mathbb{Z}
\end{equation}
Since we assumed ${\rm gcd}(n)=1$,
$C(k,n)$ is an integer.

2.
Let us write
$|n|=p^{\alpha}n'$
where $n'$ is an integer prime to $p$.
By lemma (\ref{xpb}), we have
\begin{equation}\notag
(x_1+\cdots+x_l)^{p^{i}k|n|}
=
(x_1^p+\cdots+x_l^p)^{p^{i}k|n|}
+\text{ term divisible by }p^{i+\alpha}
\qquad (i\geq 1)
\end{equation}
Then comparing the coefficient of 
$x_1^{p^{i}kn_1}\cdots x_l^{p^ikn_l}$,
we obtain
\begin{equation}\notag
\frac{(p^ik|n|)!}
{(p^ikn_1)!\cdots(p^ikn_l)!}-
\frac{(p^{i-1}k|n|)!}
{(p^{i-1}kn_1)!\cdots(p^{i-1}kn_l)!}
\equiv 0
\mod p^{i+\alpha}.
\end{equation}
Moreover, the LHS is divisible by
$|n|$ by the result of 1.
Therefore the LHS is divisible by
${\rm lcm}(|n|,p^{i+\alpha})=p^i|n|$.
\end{proof}

\begin{lemma}\label{Rt}
For $R(t)\in{\mathcal L}[t]$,
there exists $h(t)\in {\mathcal L}[t]$
such that
\begin{equation}\notag
R(t)^{p^i}=R(t_p)^{p^{i-1}}+p^ih(t).
\end{equation}
\end{lemma}
\begin{proof}
It was shown in \cite{bp}
that 
\begin{equation}\notag
t^{p^ib}=(t_p)^{p^{i-1}b}+p^i {}^{\exists}g_1(t)
\qquad(g(t)\in\mathbb{Z}[t]).
\end{equation}
Therefore, for $r(t)\in \mathbb{Z}[t]$,
\begin{equation}
r(t^p)^{p^{i-1}}-r(t_p)^{p^{i-1}}=
p^i {}^{\exists}g_2(t),
\qquad (g_2(t)\in\mathbb{Z}[t])
\end{equation}
The case $i=1$
is the above result of \cite{bp}
and
the proof is by induction on $i$.
On the other hand, 
it was shown in \cite{peng}, Lemma 5.1
that, for $r(t)\in\mathbb{Z}[t]$,
\begin{equation}\notag
r(t)^{p^i}-r(t^p)^{p^{i-1}}=
p^i {}^{\exists}g_3(t)\qquad (g_3(t)\in\mathbb{Z}[t]). 
\end{equation}
Hence, for $r(t)\in\mathbb{Z}[t]$,
\begin{equation}\label{pol}
r(t)^{p^i}-r(t_p)^{p^{i-1}}=
p^i {}^{\exists}g_4(t)\qquad (g_4(t)\in\mathbb{Z}[t]). 
\end{equation}
Let us write $R(t)=r_1(t)/r_2(t)$
where $r_1(t),r_2(t)\in\mathbb{Z}[t]$ and $r_2(t)$ is monic.
Then
\begin{equation}\notag
\begin{split}
&R(t)^{p^i}-
R(t_p)^{p^{i-1}}
=
\frac{r_1(t)^{p^i}}{r_2(t)^{p^i}}
-\frac{r_1(t_p)^{p^{i-1}}}{r_2(t_p)^{p^{i-1}}}
\\
&=
\frac{1}{r_2(t)^{p^i}r_2(t_p)^{p^{i-1}}}
\big(
r_1(t)^{p^i}(r_2(t)^{p^{i-1}}-r_2(t)^{p^i})+
r_2(t)^{p^i}(r_1(t)^{p^1}-r_1(t)^{p^{i-1}})
\big).
\end{split}
\end{equation}
By (\ref{pol}),
the numerator is written in the form
$p^i g_5(t)$ with some $g_5(t)\in\mathbb{Z}[t]$. 

The lemma is proved.
\end{proof}

\subsection{Proof of Lemma \ref{arith}}
Now we give a proof of lemma \ref{arith}.

If $k=1$, the statement is trivial,
so we prove the  $k>1$ case.
Let 
\begin{equation}\notag
k=k_1^{a_1}\cdots k_s^{a_s}
\end{equation}
be the prime decomposition of $k$.
It is sufficient to prove the following
for every $i$ ($1\leq i\leq s$).
\begin{equation}\label{whatweshow}
\begin{split}
\sum_{k';k'|k}
\mu\Big(\frac{k}{k'}\Big)
\frac{|k'{n}|!(-1)^{k'|n|}}{(k'n_1)!\cdots(k'n_l)!}
%
%
R(t_{k/k'})^{k'}
=k_i^{a_i}|n| {^{\exists}}f_i(t),
\qquad 
f_i(t)\in {\mathcal L}[t].
\end{split}
\end{equation}
%
%

Since the proof of (\ref{whatweshow}) 
is the same for any $i$ $(1\leq i\leq s)$,
we only show the case $i=1$.
Let $j=k/k_1^{a_1}$.
Note that for any divisor $k'$ of $k$,
\begin{equation}\notag
\mu\Big(\frac{k}{k'}\Big)
=\begin{cases}
\mu\big(j/j'\big)& k'=k_1^{a_1}j'\\
-\mu\big(j/j'\big)& k'=k_1^{a_1-1}j'\\
0&\text{ otherwise}
\end{cases}
\end{equation}
Therefore
\begin{equation}\notag
\begin{split}
&\text{LHS of (\ref{whatweshow})}
\\&=
\sum_{j';j'|j}
\mu\Big(\frac{j}{j'}\Big)
\Bigg[
\frac{|{k_1^{a_1}j'n}|!(-1)^{k_1^{a_1}j'|n|}}
{(k_1^{a_1}j'n_1)!\cdots(k_1^{a_1}j'n_l)!}
R(t_{j/j'})^{k_1^{a_1}j'}
-
\frac{|{k_1^{a_1-1}j'n}|!(-1)^{k_1^{a_1-1}j'|n|}}
{(k_1^{a_1-1}j'n_1)!\cdots(k_1^{a_1-1}j'n_l)!}
R(t_{k_1 j/j'})^{k_1^{a_1-1}j'}
\Bigg]
\\
&=\sum_{j';j'|j}
\mu\Big(\frac{j}{j'}\Big)
(-1)^{k_1^{a_1}j'|n|}
R(t_{j/j'})^{k_1^{a_1}j'}
\underbrace{
\Bigg[
\frac{|{k_1^{a_1}j'n}|!}
{(k_1^{a_1}j'n_1)!\cdots(k_1^{a_1}j'n_l)!}
-
\frac{|{k_1^{a_1-1}j'n}|!}
{(k_1^{a_1-1}j'n_1)!\cdots(k_1^{a_1-1}j'n_l)!}
\Bigg]
}_{\star}
\\
&+\sum_{j';j'|j}
\mu\Big(\frac{j}{j'}\Big)
\underbrace{
\frac{|{k_1^{a_1-1}j'n}|!}
{(k_1^{a_1-1}j'n_1)!\cdots(k_1^{a_1-1}j'n_l)!}
}_{\ast}
\underbrace{
\Bigg[
\Big((-1)^{j'|n|}
R(t_{j/j'})^{j'}\Big)^{k_1^{a_1}}
-\Big((-1)^{j'|n|}
R(t_{k_1j/j'})^{j'}\Big)^{k_1^{a_1-1}}
\Bigg]
}_{\dagger}
.
\end{split}
\end{equation}
By lemma \ref{art}, $(\star)$ is divisible by $k_1^{a_1}|{n}|$
and
($\ast$) is divisible by $|{n}|$.
By lemma \ref{Rt},
there exists $h(t_{j/j'})\in {\mathcal L}[t_{j/j'}]$
such that 
($\dagger$) is written as
$k_1^{a_1} h(t_{j/j'})$.
Since $ {\mathcal L}[t_{j/j'}]\subset {\mathcal L}[t]$,
$h(t_{j/j'})\in {\mathcal L}[t]$.
Therefore (\ref{whatweshow}) is proved.

This completes the proof of lemma \ref{arith}.
\section{Multiplicativity of $\Psi$}
\label{appendix:Psi}
In this section, we show that
the map $\Psi$, introduced in the proof of
proposition \ref{free},
is multiplicative with respect to the graph union.

Since the combined amplitudes of two equivalent VEV
forests are the same,
(\ref{def:Psi}) is equivalent to:
\begin{equation}\label{PsiG}
\Psi(G)=\#\aut(G)\cdot N(W)
\frac{1}
{z_{\Vec{\mu}}z_{\Vec{\nu}}}
{\mathcal H}(W)\Vec{Q}^{\Vec{d}}
\qquad 
(W \text{ such that } \overline{W}=G),
\end{equation}
where
$N(W)$ is
the number of combined forests equivalent to
$W$.

$N(W)$ is described as follows.
It
is equal to the number of ways to assign
leaf indices to leaves of
$G\in \overline{{\rm Comb}}_{\Gamma}^{\bullet}(\vec{\mu},\vec{\nu})$.
Therefore
\begin{equation}\label{Nw}
N(W)=N_0(G)
\prod_{T:\text{VEV tree in $W$}}
N'(T),
\end{equation}
where $N'(T)$
is the number of
VEV trees
equivalent to $T$
and
$N_0(G)$
is the number of ways to distribute the 
leaf indices to 
VEV trees 
in $G$.
%
The latter is
\begin{equation}\label{Nzero}
\begin{split}
N_0(G)&=\frac{\#\aut(\Vec{\mu})\,\#\aut(\Vec{\nu})}
{\#\aut(G)}
\prod_{T:\text{VEV tree in $G$}}
\Big(
\prod_{f\in F_3(\Gamma)\setminus \hat{F}_3(\Gamma)}
\#\aut(\nu^f(T))
\Big)^{-1}
\\&\times
\prod_{\begin{subarray}{c}\{T,T'\};\\
T,T' :\text{VEV trees in $G$},\\
T\neq T'
\end{subarray}}
\Big(
\prod_{v\in V_3(\Gamma)}\#\aut(\mu^v(T,T')
\prod_{f\in \hat{F_3}(\Gamma)}\#\aut(\nu^f(T,T'))
\Big)^{-1},
\end{split}
\end{equation}
where
\begin{equation}\notag
\begin{split}
\nu^f(T)
&=\{\nu^f_i|\text{ $\nu^f_i$ is associated to
a leaf of $T$}\}
\qquad
(f\in F_3(\Gamma)\setminus \hat{F_3}(\Gamma)),
\\ %
\mu^v(T,T')
&=\{\mu^v_i|\text{ $\mu^v_i$ is associated to 
a bridge joining $T$ and $T$'}\} 
\qquad(v\in V_3(\Gamma)),
\\
\nu^f(T,T')  
&=\{\nu^f_i|\text{ $\nu^f_i$ is associated to a bridge joining
$T$ and $T'$}\}
\qquad (f\in \hat{F}_3(\Gamma)).
\end{split}
\end{equation}
Substituting (\ref{Nw})(\ref{Nzero}) into (\ref{PsiG}),
we find out that
$\Psi$ is multiplicative.
%
%

\end{document}